\documentclass[review,3p]{elsarticle}

\usepackage{lineno,hyperref}
\modulolinenumbers[5]

\biboptions{authoryear}

\usepackage{ulem,amsmath, amsfonts,amssymb,graphicx, color, rotating, url, tikz, pgffor,xkeyval, listings, mathtools, etoolbox, bm, aliascnt, epstopdf, multirow, hyperref, cleveref, subcaption, float, pgfgantt, pgfplots, booktabs, mathrsfs, palatino, multicol, gensymb}
\usetikzlibrary{arrows, automata, shapes, decorations, shapes.geometric, positioning, patterns}
\pgfplotsset{compat=1.14}

\usepackage{algorithm}
\usepackage{algorithmic} 
\usepackage{algorithmic,eqparbox,array}

\usepackage[colorinlistoftodos,prependcaption,textsize=scriptsize]{todonotes}

\usepackage[nolist]{acronym}
\begin{acronym}
\acro{SFSMP}{Stochastic Fleet Size and Mix Problem}
\acro{TCO}{Total Cost of Ownership}
\acro{SESO}{Simulation-Enabled Strategic Optimizer}
\acro{SAA}{Sample Average Approximation}
\acro{ALNS}{Adaptive Large Neighborhood Search}
\acro{ICEV}{Internal Combustion Engine Vehicle}
\acro{EV}{Battery Electric Vehicle}
\acro{FSMVRP}{Fleet Size and Mix Vehicle Routing Problem}
\acro{SOC}{State of Charge}
\acro{SISR}{Slack Induced String Removal}
\end{acronym}

\usepackage{setspace}
\begin{document}
\DeclarePairedDelimiterX{\abs}[1]{\lvert}{\rvert}{\ifblank{#1}{{}\cdot{}}{#1}}
\DeclarePairedDelimiter{\ceil}{\lceil}{\rceil}
\DeclarePairedDelimiter{\floor}{\lfloor}{\rfloor}
\begin{frontmatter}

\title{Stochastic Fleet Mix Optimization: Evaluating Electromobility in Urban Logistics}


\author[mymainaddress,mysecondaryaddress]{Satya S.~Malladi\corref{mycorrespondingauthor}}
\cortext[mycorrespondingauthor]{Corresponding author}
\ead{samal@dtu.dk}

\author[mymainaddress]{Jonas M.~Christensen} \ead{jomc@dtu.dk}
\author{David Ramírez}

\author[mymainaddress]{Allan Larsen}\ead{alar@dtu.dk}
\author[mymainaddress]{Dario Pacino}\ead{darpa@dtu.dk}

\address[mymainaddress]{DTU Management, Technical University of Denmark,  Denmark}

\begin{abstract}
In this paper, we study the problem of optimizing the size and mix of a mixed fleet of electric and conventional vehicles owned by firms providing urban freight logistics services. Uncertain customer requests are considered at the strategic planning stage. These requests are revealed before operations commence in each operational period. At the operational level, a new model for vehicle power consumption is suggested. In addition to mechanical power consumption, this model accounts for cabin climate control power, which is dependent on ambient temperature, and auxiliary power, which accounts for energy drawn by external devices. We formulate the problem of stochastic fleet size and mix optimization as a two-stage stochastic program and propose a sample average approximation based heuristic method to solve it. For each operational period, an adaptive large neighborhood search algorithm is used to determine the operational decisions and associated costs. The applicability of the approach is demonstrated through two case studies within urban logistics services.
\end{abstract}

\end{frontmatter}


\noindent\textbf{Funding: }The funding body will be acknowledged following peer review.
\section{Introduction}
Urban freight transport logistics is an energy-intensive activity that involves vehicle movement on congested roadways in densely populated regions. The widespread use of conventionally powered \acp{ICEV} over many years has been negatively impacting the environment in many ways. Greenhouse gas emissions from \acp{ICEV} not only aggravate climate change but also take a serious toll on the cardio-pulmonary and respiratory health of the humans inhaling them \cite{Heinrich2004TrafficDisease}. In addition, prolonged exposure to high noise levels due to the operation of \acp{ICEV}
is detrimental to health \cite{Kijewska2016FreightExample}. 
With an increasing proportion of population living in urban areas (predicted to grow from the current 55\% to 68\% in 2050 \citep{UN201868Affairs}) and burgeoning demand for freight transport services in cities, owing largely to the spurt of e-commerce, the adversity resulting from the use of \acp{ICEV} will only worsen the livability of urban spaces. It is now an opportune moment for all stakeholders involved in various degrees of decision-making to seek win-win measures to ensure sustainability \cite{CIVITAS2015MakingSustainable}.
Given the short driving distances in urban localities, using \acp{EV} for goods distribution in urban areas is a pertinent approach available to firms for moving towards sustainable transport operations. The electrification of urban freight transportation fleet is a high impact solution in the context of fighting climate change. 
While \acp{EV} may still indirectly produce carbon dioxide emissions due to the widespread use of fossil fuel based electricity, they diminish urban particulate as well as noise pollution significantly. Ambitious environmental regulations, together with the technological advances made by the industry in recent years, have made electro-mobility a real alternative for companies providing services in urban localities.

Although recent major improvements in Lithium-ion batteries, both in performance and production costs, have heralded accelerated integration and improved market share of \acp{EV} \cite{Waraich2015, InternationalRenewableEnergyAgency2017ElectricBrief}, firms are weary of the disruption of operations due to a problematic mismatch between planned driving range and realized driving range of \acp{EV} during fleet operations. 
The variation of ambient temperature during the planning horizon and the difficulty of estimating the auxiliary energy usage at the planning stage are major factors contributing to this mismatch, signaling the need for more advanced models for estimating energy consumption in the fleet electrification research literature.

\begin{figure}[ht!]
\begin{center}
\includegraphics[scale = 1.6]{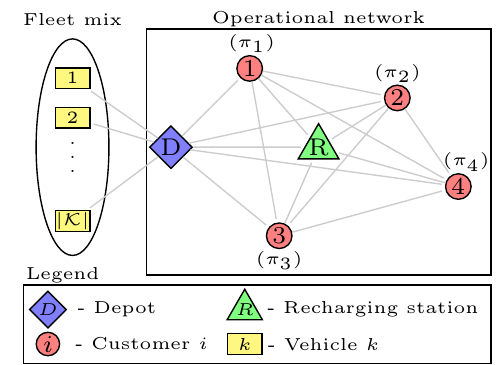}
\caption{An example network for the \ac{SFSMP}}
\label{sfsmp_network}
\end{center}
\end{figure}

In this paper, we investigate the smart adoption, integration, and the efficient use of electro-mobility in urban logistics. In particular, we study the strategic problem of identifying the size of a fleet of mixed \acp{EV} and/or \acp{ICEV} for companies involved in urban logistics. We assume customer requests and temperature are stochastic at the strategic planning stage and are only revealed prior to operations everyday. The operational problem is a vehicle routing problem with the realized requests. Formally, we refer to this problem as the \ac{SFSMP}. We note that the quantity of demand at each customer is known with certainty given the customer presents a request.
\Cref{sfsmp_network} gives a visual presentation of the \ac{SFSMP}, where customers, charging stations, a depot, and a vehicle fleet are depicted. The nodes representing the fleet mix show the complete set of vehicles considered. At the strategic level, the optimization decisions concern the selection of which vehicles should compose the fleet. The decision is represented by the edges connecting each vehicle with the depot. The network, between the customers and the depot, represents the operational periods where customer requests are stochastic and appear with probability $\pi_i$. The operational period is also affected by another stochastic variable, the ambient temperature.
We argue that determining the optimal solution of the \ac{SFSMP} is currently intractable. We thus propose a heuristic approach based on \acf{SAA} \cite{Ahmed2002TheRecourse} to solve the \ac{SFSMP} and thereby determine the fleet requirement over the entire horizon. 

There has been a significant quantum of literature on the fleet size and mix problem with electric vehicles when demands and requests are deterministic and constant over the entire planning horizon \cite{Golden1984TheProblem,Hiermann2016}. However, to the best of our knowledge, the fleet size and mix literature on urban freight logistics has not considered the uncertainty of operational demand requests in strategic planning. 

The contributions of our work are several. First, we propose the \ac{SFSMP}, a novel stochastic optimization problem for strategic fleet size and mix decision-making  in urban freight logistics. Second, we upgrade the state-of-the-art energy consumption model presented by \citet{Goeke2015RoutingVehicles} by explicitly accounting for the energy required for auxiliary usage (such as heating and cooling the cabin space of vehicles,  connecting external devices, etc.). Third, we propose a generalization of the \ac{FSMVRP} for the evaluation of the operational problem appearing in the \ac{SFSMP} that can handle time windows of service, compatibility constraints between vehicles/drivers and customers, driving range limitations, and en-route charging necessities. Fourth, we propose a heuristic variation of the \ac{SAA} method to solve the \ac{SFSMP}, where the decision-making in each operational period (both within the strategic decision-making process as well as outside) is carried out using an adaptation of the state-of-the-art \acf{ALNS} algorithm. Finally, we demonstrate the applicability of the approach on two different case-studies of service operations in urban freight logistics.

The remainder of this paper is organized as follows. In  \Cref{lit_rev}, we present a summary of the literature on fleet sizing problems, urban freight logistics with \acp{EV}, energy consumption models, and discuss solution methods in the context of freight logistics planning. A mathematical formulation is presented for the strategic optimization problem with an embedded vehicle routing model for the operational optimization problem in  \Cref{prob_form}. 
We present our adaptation of \acs{SAA} for solving the \ac{SFSMP}  in \Cref{methodology}. In \Cref{cases}, we present case studies of two Danish organizations embracing electric mobility through research initiatives. We summarize our insights from the analysis of the novel \ac{SFSMP} and the practical cases in \Cref{conclusion}.     
 
 \section{Literature Review \label{lit_rev}}
Our review focuses on the four major aspects of this paper, namely, fleet size and mix optimization, commercial use of electric vehicles in urban logistics, energy consumption modeling, and relevant methods for solving the problem of interest. 
 \subsection{Fleet size and mix optimization}
The fleet size and mix vehicle routing problem (FSMVRP) with deterministic demands that are identical on all operational days accounting for fixed leasing costs and variable routing costs  was formalized by \citet{Golden1984TheProblem}. It was solved using savings based constructive heuristics and a two-phase route first and cluster second algorithm \cite{Golden1984TheProblem}, a generalized assignment based heuristic \cite{Gheysens1984AProblems}, a matching based tour - fusion heuristic \cite{Desrochers1991AProblem}, tabu-search based local search methods \cite{Taillard1999AVRP,Gendreau1999AProblem,Osman1996LocalProblem}, and a sweep-based heuristic method that can also solve non-Euclidean problems \cite{Renaud2002AProblem}. \citet{Baldacci2009ValidCosts} proposed a two-commodity network formulation of the FSMVRP and presented various valid inequalities to strengthen the formulation.

\citet{Liu1999TheWindows} extended this problem to include hard time windows (FSMVRPTW). The FSMVRPTW was solved in the literature using an insertion-based heuristic and an improvement scheme \cite{Liu1999TheWindows}, a sequential insertion heuristic \cite{Dullaert2002NewWindows},  and a ruin and recreate method \cite{DellAmico2007HeuristicWindows}. More recently, 
\citet{Hiermann2016} generalized this problem to incorporate driving range constraints of electric vehicles and the option of en-route recharging at the operational level. They solved the problem using an exact branch and price procedure and a hybrid \ac{ALNS} based heuristic demonstrating the effectiveness of their method on a new set of  benchmark instances. The approach of simulation has been used in \cite{Cortes2011ADemand} to sample various scenarios to aid in the daily dispatch of technicians.

 All the above fleet size and mix problems studied in the literature have considered deterministic customer sets that are known at the strategic planning stage and remain identical through out the strategic planning  horizon. \citet{Pasha2016SimpleProblem} studied a FSMVRP (without time windows) with deterministic demands that vary on different days.  A heuristic method was used to obtain a good common fleet size and mix across the all days after determining the best fleet size and mix independently for each day. Fleet size and mix problems appearing in maritime logistics have considered uncertainty in two-stage stochastic programming problems, such as maintenance operations planning in offshore wind farm management \cite{Stalhane2019OptimizingFarms} and fleet renewal planning \cite{Bakkehaug2014AShipping} under uncertain weather conditions. However, uncertain customer requests at the strategic planning stage have not been considered in the fleet size and mix literature on urban freight logistics so far. In the current paper, we study a fleet size and mix problem in which customer requests are not revealed until the operational period begins.  

 Fleet sizing boils down to a tradeoff between acquisition costs and operational costs, thus rendering an accurate consideration of the \ac{TCO}  crucial. The \ac{TCO} of a fleet of vehicles owned by a company must account for various fixed costs (such as capital cost, insurance, and taxes) and variable costs that depend on the duration of operation (such as maintenance and repair costs and energy consumption cost).  \citet{Rogge2018ElectricInfrastructure} and \citet{Schiffer2018ElectricNetworks} performed \ac{TCO} analysis for commercial mixed electric fleets in the context of  public transit and mid-haul logistics respectively.    
 In the current paper, we consider a \ac{TCO} that comprises of one-time purchase costs at the strategic level and energy and maintenance costs that are proportional to the travel time in addition to penalty for not serving some customer in every operational period. A more detailed description of the \ac{TCO} can also be handled by our mathematical model.   
 \subsection{Urban logistics with electric vehicles}
 \citet{Pelletier2016} provided an in-depth review of electric vehicle technologies and summarized  transportation science literature on goods distribution with \acp{EV}.  \citet{Crainic2009ModelsSystems} presented a brief survey of city logistics problems and discussed  models and solution methods for two-tiered city logistics systems.  \citet{Savelsbergh201650thOpportunities} delivered a contemporary perspective of the changing face of city logistics due to the advent of e-commerce prevalence, customer impatience, collaborative consumption, and sustainability consciousness paired with the advancement of digital as well as physical technologies.   \citet{Tang2019TheEra} analyze the importance of logistics function in the fourth industrial revolution and offer insights into how the logistics industry can generate economic, environmental, and social value by transforming itself suitably. A common thread revealed in these papers is the need to manage electric vehicles in urban freight logistics.
Tracing the beginning of electric VRPs, \citet{Conrad2011TheProblem} proposed the recharging VRP and were followed by \citet{Erdogan2012AProblem} with a green vehicle routing problem, and by \citet{Schneider2014TheStations} with the EVRP with time windows and recharging stations. The literature in this area has come a long way, with a more recent variant by \citet{Cortes-Murcia2019TheCustomers} allowing alternate modes of transport for customer visits while recharging to reduce loss of productivity.
 \citet{Lin2016} formulated an electric vehicle routing problem with a heterogeneous fleet that accounts for the impact of speed and vehicle load on battery consumption, allowing en-route charging. 
  \citet{Goeke2015RoutingVehicles} generalized \citep{Schneider2014TheStations} by considering non-linear expressions for fuel and battery consumption and planning with a mixed fleet of \acp{EV} and \acp{ICEV}. \citet{Hiermann2019RoutingVehicles} proposed an integrated EVRPTW with a mixed fleet consisting of \acp{EV}, plug-in hybrid vehicles, and \acp{ICEV}, using constant rates of battery and fuel consumption in their model. Decision-making for green transportation has, in the last decade, helped reducing the negative externalities of freight transport (see \cite{Bektas2016GreenRouting, Bektas2011TheProblem, Bektas2018TheTransportation}). We note that the problem of electrification for urban service logistics  is significantly relevant given the prevalent interest on adoption of \acp{EV}.
 \subsection{Energy consumption models}
 Previous works by \citet{Erdogan2012AProblem} and \citet{Schneider2014TheStations}  considered linear energy consumption functions that were linear in the distance travelled and speed. More advanced energy accounting was performed for \acp{ICEV} in  \citet{Demir2012AnProblem}. In their paper, mechanical energy consumption was characterized to overcome rolling resistance, aerodynamic resistance, and gravitational force using expressions that are non-linear in speed. \citet{Goeke2015RoutingVehicles} extended this model to model battery energy consumption of \acp{EV} in a load-dependent fashion. \citet{Murakami2017ARouting} use this model to present a vehicle routing problem with electric and diesel vehicles with the notion of an original graph for estimating at more accurate road slopes, wait times at road intersections, etc. 
 The existing energy consumption models do not consider the energy consumption due to cabin heating and cooling, which contribute significantly to the reduced driving range during peak winter and peak summer. The power consumed for cabin heating at 0\degree C and that for cabin cooling at 30\degree C is around 40\% of the total vehicle power for a mid-sized pickup van (seen later in \Cref{temp_power_kangoo}). In the current paper, we extend the model in \citep{Goeke2015RoutingVehicles} further to include auxiliary energy consumption due to heating and cooling loads and use the non-linear recharging function proposed by \citet{Montoya2017} to determine the time required to charge a given amount of energy. Another aspect in energy consumption calculation of \acp{EV} is the effect of battery degradation \cite{Pelletier2017BatteryModels}. We ignore its impact in the current paper.

 \subsection{Solving the \acf{SFSMP}}
 The literature on planning problems concerning the strategic and tactical decisions that account for operational efficiency is sparse. \citet{Dempster1981AnalyticalSystems} presented a framework for modeling hierarchical planning problems as multi-stage stochastic programs that are capable of capturing uncertainty at the lower levels of decision-making. More recently,  \citet{Crainic2015ModelingPlanning} considered the tactical planning of the design of the service network in the first tier while accounting for uncertainty in demands in the second tier in a two tier city logistics system. The performance of four tactical plans with associated recourse strategies in the operational stage was evaluated in their paper. 
 To solve two-stage stochastic programs with integer recourse, seminal works have appeared in the literature bringing to the fore various approaches \cite{Shapiro20009LecturesProgramming}  including \acs{SAA} \cite{Ahmed2002TheRecourse},  which was proved to result in a solution that converges to the solution of the original problem as the sample size increases \cite{Schultz1996RatesRecourse}. In the current paper, we employ \acs{SAA} to solve the strategic fleet size and mix problem 
 using the metaheuristic  \ac{ALNS} to solve the operational problem. 
 The metaheuristic method \ac{ALNS} was first introduced by \citet{Pisinger2007AProblems} and is considered to be one of the most efficient methods for solving vehicle routing problems with deterministic requests. Other metaheuristics include a genetic algorithm based method \cite{Prins2004AProblem} and tabu search methods \cite{Gendreau2008AProblem}.

 In the next section, we present the formulation of the strategic problem as a two-stage stochastic program. 
\section{Problem Formulation\label{prob_form}}

 
 In this paper, we address the introduction and integration of \acp{EV} in a fleet of commercial vehicles operated within urban areas. We consider the vehicles to be used for either cargo pickup and/or service actions (e.g. electric installations, medical visits etc.). The \ac{SFSMP} is defined as a strategic decision  problem. It aims at minimizing the \ac{TCO} while determining the optimal fleet size and mix for service operations repeating over multiple operational periods. We assume \textit{a posteriori} customer requests that are not known, owing to the lack of advance information about each day's requests at the strategic level. However, at the operational level, requests are deterministic. Since strategic decision-making here relies on the trade-off between fleet acquisition cost and average operational cost, an accurate estimation of the operational costs is needed as well.  

We formulate the  \ac{SFSMP} as a two-stage stochastic mixed integer program. The strategic horizon consists of a number of operational periods. For example, there may be 1200 days of operation in a strategic horizon of 5 years. 
The first stage problem is the {strategic problem} of determining the optimal fleet size and mix that minimizes the sum of the fleet acquisition cost and the expected second stage cost. The strategic decision of the fleet mix is determining the set of vehicles acquired to serve the requirements of the urban logistics operations over the entire horizon in the form of $\bm{w}$, a vector of binary variables such that $w_k$ takes a non-zero value if vehicle $k$ of the master list of vehicles $\mathscr{K} $ is acquired.  The second stage problem consists of multiple operational periods in which routing decisions of each period are affected by the first stage fleet mix decision only and not the routing decisions of the other periods. In each period of the second stage, only the vehicles chosen for acquisition in the first stage can be used to serve realized customer requests. In every operational period, the {operational problem} must be solved, which consists in the routing of  the drivers, and can be formally defined as a \textit{mixed electric fleet vehicle routing problem with time windows and compatibility constraints}. We consider a heterogeneous fleet of vehicles of different types varying in size and/or energy source. Although the feasible driving range of most commercial \acp{EV} available in the market is typically higher than the distance traversed in urban logistics tours, we allow the possibility of making en-route recharging visits. Additionally, since the duration available for providing service in urban logistics during the day is limited to a daily work shift in a business to business (B2B) service environment, we limit the number of potential visits to recharging stations to one visit. The set of all recharging stations is given by $F$. We remark that adding recharging stops to the routes would result in a loss of productivity even with the latest charging technology in the context of urban logistics and tight time constraints. 
Each customer request is characterized by a known demand   quantity, a hard time window of service, and is associated with a specific level of skill requirement. Without loss of generality, we focus on pickup demands in the case studies that follow. Each driver is assigned a specific skill set that must include the customer's skill requirement for service to take place. Therefore,  we model compatibility constraints between drivers and customer requests. We note that this problem is generic and is capable of handling special cases where all drivers are compatible with all requests. 

\subsection{First stage problem}
 The formulation of the first stage problem of the \ac{SFSMP} is presented with the definition of the sets and parameters  presented in  \Cref{first_stage}. In this stage, the decision variables are binary to indicate whether or not vehicles in the master list should be acquired. Bold notation is used for vector quantities. The vector $\bm{w}$ indicates whether or not each vehicle in the master list of vehicles $\mathscr{K} $  is acquired.

\begin{table}[ht!]
\flushleft
\caption{Sets and parameters used in the formulation of the first stage problem \label{first_stage}}
\begin{tabular}{l|l}
\toprule 
\multicolumn{2}{c}{\textbf{Sets}}\\
 $\mathscr{K}$ & Master set of all possible vehicles    \\
$ \mathcal{N}$ & Master list of all customers \\ 
$r_i$ & Bernoulli random variable for the occurrence of request from  customer  $i$ with state space $\{0,1\}$\\ 
$\lambda_i$ & a realization of request random variable $r_i$\\ 
 $\pi_i$&Probability of the occurrence of request from customer $i$  \\
$\Omega(\mathcal{T})$ & State space of the random  ambient temperature  \\ 
 \multicolumn{2}{c}{\textbf{Parameters}}\\
$\beta_k$ & Acquisition cost of vehicle $k$ \\ 
$n_o$& Number of varieties of operational periods    \\ 
$n'$& Number of operational period of each variety    \\ 
$ \mathcal{T}$ & Ambient temperature; a random variable in the first stage \\ 
$T$ & A realization of the ambient temperature \\
\bottomrule
\end{tabular}
\end{table}

The formulation of the \ac{SFSMP} is presented below. The strategic objective function accounts for a) the cost of acquiring the vehicles  constituting fleet mix configuration determined by the decision vector $\bm{w}$ and b) the total expected operational cost of serving appeared requests with the acquired fleet. The set of acquired vehicles determined by $\bm{w}$ is a subset of the master set of vehicles $\mathcal{K}$.  In the expression below, the first term captures the fleet acquisition cost using vehicle acquisition costs $\beta_k$ while the second term captures the total expected operational cost summed over $n_o$ varieties of operational periods (such as summer, winter, etc.) with $n'$ operational periods in each variety. In the second term, $\mathcal{C}(\bm{w}, \bm{\lambda}, T)$ represents the deterministic operational cost  of serving the set of realized customers determined by the vector of random variables $\bm{r}$ (or binary vector $\bm{\lambda}$ of corresponding realized values) at temperature $\mathcal{T}$ which is a random variable (or realized temperature $T$) with fleet determined by$\bm{w}$. Each $\bm{r}^j, j \in \{1,\dots, n_o\}$ represents the binary request vector of customers in the $j$th variety of operational periods. Similarly, $\Omega(\mathcal{T})$ represents the state space of the ambient temperature $\mathcal{T}$ which is a random variable. 
\begin{equation} \min_{\bm{w} } \bigg\{ \sum\limits_{k \in \mathscr{K}} \beta_k w_k + \sum\limits_{\bm{r} \in \{\bm{r}^1 , \dots, \bm{r}^{n_o} \}} n' \ \mathbb{E}_{\bm{r}, \mathcal{T}} \left [ \mathcal{C} (\bm{w}, \bm{r}, \mathcal{T})  \right ]  \bigg \},  \label{strat_obj_0} 
\end{equation}
If it is  assumed that each operational period faces identical uncertainty in the occurrence of customer requests  and ambient temperature, the strategic objective is re-written as \eqref{strat_obj}. 
\begin{equation} \min_{\bm{w}} \bigg\{ \sum\limits_{k \in \mathscr{K}} \beta_k w_k + n_o n' \  \mathbb{E}_{\bm{r}, \mathcal{T}} \left [ \mathcal{C}(\bm{w}, \bm{r}, \mathcal{T})  \right ]  \bigg \},  \label{strat_obj} 
\end{equation}
where  \begin{equation} \mathbb{E}_{\mathcal{N}, \mathcal{T}} \left [ \mathcal{C}(\bm{w},  \bm{r}, T)  \right ]  = \sum_{\bm{\lambda} \in \{0,1\}^{ \mid\mathcal{N}\mid}} \sum_{ T \in \Omega(\mathcal{T}) } \text{Pr}(\bm{\lambda},T) \  \mathcal{C}(\bm{w}, \bm{\lambda}, T) \end{equation} 
and $\text{Pr} (\bm{\lambda},T) = \text{Pr} (T) \ \prod_{i\in \mathcal{N}} \pi_i^{\lambda_i}(1-\pi_i)^{1 - \lambda_i}$,   assuming mutual independence between the occurrence of customer requests and the temperature in the operational period, with the probability of the occurrence of request $i$ given by $\pi_i$. In every operational period, a rich vehicle routing problem with deterministic request set determined by $\bm{\lambda}$ and ambient temperature $T$, formulated in the following subsection, must be solved, given that the set of acquired vehicles determined by $\bm{\omega}$. We note that the number of such request sets is exponential in the cardinality of  $\Omega(\mathcal{N})$.   
\subsection{Second stage problem}
\begin{table}[ht!]
\flushleft
\caption{Sets and parameters used in the formulation of the second stage problem \label{second_stage}}
\begin{tabular}{l|l}
\toprule
 \multicolumn{2}{c}{\textbf{Sets}} \\
 \\ \midrule
$A$ & Set of arcs $(i,j), i,j \in N'_{0,n+1}$ in the graph \\
$F$ & Set of recharging stations  \\
$ \bm{\omega}$ & A known vector of membership of the master set of vehicles $\mathscr{K}$ \\
 $K$ & The set of vehicles used in the second stage problem as determined by $ \bm{\omega}$  \\
 $K^{E} \subseteq K$& Subset of vehicles of $K$ that are \acp{EV} \\
 $\bm{\lambda}$ & A binary vector of realized request random variables  \\
 $N$& The set of customers $\{1,2,...,n\}$ determined by $\bm{\lambda} $; a subset of $\mathcal{N} $ \\
 $N'$ & $N \cup F$ \\
  $N'_{0}$& $N\cup F\cup \{0\} $\\
 $N'_{n+1}$& $N\cup F\cup \{n+1\} $\\
 $N'_{0,n+1}$& $N\cup F\cup \{0,n+1\} $ \\
 \midrule
 \multicolumn{2}{c}{\textbf{Parameters}} \\
\midrule
$0, n+1$& Depot vertices \\
 $a_{ik}$& 1 if task to be performed at vertex $i$ is compatible with vehicle $k$, 0 otherwise \\
$b_i (y)$ & Function to determine the time taken to recharge $y$ units of charging at station $i$\\
 $c_k$& Energy cost of vehicle $k$ in USD/kWh \\
  $d_i$& Demand of customer $i$ given the request appears \\
 $[e_i,l_i]$& Time window of vertex $i$ \\
 $E_k$& Battery range of vehicle $k$ in kWh \\
 $f_k$& Daily maintenance and usage cost of vehicle $k$ \\
 $g_i$& Penalty or taxi cost for not serving customer $i$  in USD \\
 $M_k$ & maintenance cost  in USD per km of vehicle $k$ \\
 $p_{ijk}(q,T)$& Power consumption (kW)  on arc ($i$,$j$) for vehicle $k$ with load $q$ at temperature $T$ \\
  $Q_k$& Capacity of vehicle $k$ \\
 $s_i$& Service time at vertex $i$ \\
 $t_{ijk}$& Travel time in hours between vertices $i$ and $j$ for vehicle $k$ \\
  $v_k$& Average operational speed of vehicle $k$ in kmph \\
  \bottomrule
\end{tabular}
\end{table}
To complete the model of the \ac{SFSMP}, the deterministic operational cost $\mathcal{C}(\bm{\omega}, \bm{\lambda}, T)$ of serving customer set determined by $\bm{\lambda}$ at temperature $T$ with fleet determined by fleet membership vector $ \bm{\omega}$ will be defined for all $ \bm{\omega} \in \{0,1\}^{\mid \mathscr{K}\mid}, \ \forall \ \bm{\lambda} \in \{0,1 \} ^{\mid \mathcal{N} \mid}, \  T \in \Omega(\mathcal{T}) $. 
For ease of exposition, we let $K$ be the subset of vehicles in the master list $ \mathscr{K}$ determined by a fleet membership vector $\bm{\omega}$ and we let $N$ be the subset of customers determined by a realized request vector $\bm{\lambda}$.
Since the number of available vehicles in each period is fixed, we allow for unserved customers and impose a cost associated with them  \cite{Lau2003VehicleVehicles}. The notation used in this formulation of the operational problem is influenced by that used in the model of \citet{Goeke2015RoutingVehicles} and is presented in \Cref{second_stage} and \Cref{dec_var}. In this formulation, a network with the set of nodes $N_{0,n+1}'$ and the set of arcs $A$ is modeled. The set of nodes includes the set of customers $N$, the set of recharging stations $F$, the source depot node $0$, and the sink depot node $n+1$. A set of acquired vehicles $K$ that contains a set of \acp{EV} $K^E$ is considered given. For every node $i$, compatibility with drivers $a_{ik}$,  demand $d_i$ (set to zero at non-customers),  time window of service $[e_i, l_i]$, penalty $g_i$ for not serving it, and duration of service $s_i$ are specified. For every vehicle $k$, driving range in kWh $E_k$ (for \acp{ICEV} also), daily maintenance and usage cost $f_k$, maintenance cost per km $M_k$, load capacity $Q_k$, average operational speed $v_k$, and cost of energy in USD/kWh $c_k$ are given. Additionally, $b_i(y)$ is a function that calculates the amount of time taken to recharge $y$ units at station $i$. On every arc $(i,j)$, travel time and load-and-temperature-dependent power consumption of vehicle $k$ are given by $t_{ijk}$ and $p_{ijk}(q,T)$ respectively. 

Any binary decision variable  $x_{ijk} $ takes the value 1 if arc $(i,j)$ is used by vehicle $k$. The binary decision variables $z_k$ and $u_i$ assume the value 1 if vehicle $k$ is used and if customer $i$ is not served respectively. The continuous variable $y_{ik}$  represents the \acf{SOC} of vehicle $k$ upon entering any node $i$ while the variable $Y_{ik}$ is the \ac{SOC} of vehicle  $k$ upon leaving recharging station $i$. The arrival time of vehicle $k$ at each node $i$ is given by $\tau_{ik} $ and the remaining capacity available in vehicle $k$  upon arriving at node $i$ is given by $q_{ik}$.

 \begin{table}[ht!]
\flushleft
\caption{Decision variables in the model of the second stage problem \label{dec_var}}
\begin{tabular}{l|l}
\toprule
\textbf{Variable} & \textbf{Description} \\
\midrule
 $q_{ik}$& Remaining cargo level on arrival at vertex $i$ on vehicle $k$ \\
 $\tau_{ik}$& Arrival time at vertex $i$ \\
 $u_i$& 1 if customer $i$ is not served, 0 otherwise \\
 $x_{ijk}$& 1 if arc $(i,j)$ is traversed by vehicle $k$, 0 otherwise \\
 $y_{ik}$& Remaining battery capacity on arrival at vertex $i$ by vehicle $k$ \\
 $Y_{ik}$& Remaining battery capacity on departure from recharging station $i$ by vehicle $k$ \\
 $z_k$& 1 if vehicle $k$ is used, 0 otherwise \\
  \bottomrule
\end{tabular}
\end{table}

 The objective function \eqref{objective} minimizes the total fleet operation cost, which is composed of variable
 energy and maintenance costs, fixed usage costs that can account for acquisition, insurance, maintenance, and driver salaries, and  penalties for not serving customers.
 
 Constraints \eqref{getunserved} determine the customers that are not served and ensure that all customer nodes are visited at most once. Constraints \eqref{atmost1} allow at most one recharging visit on every electric vehicle tour. We add this restriction since the productivity of vehicles and drivers performing intra-day commercial services with short driving duration will take a hit when multiple recharging visits take place during the day. Constraints \eqref{flowEV} and \eqref{flowICEV} ensure flow conservation at the non-depot vertices. Constraints \eqref{NoFvisitICEV} ensures that \acp{ICEV} cannot travel to or from recharging stations. This constraint is not included in the formulation of \citet{Goeke2015RoutingVehicles}. We present examples that illustrate the need for this constraint in \Cref{App:ConstraintExplanation}. Constraints \eqref{getvehicles} determine the vehicles used for operations. Constraints \eqref{comp} are compatibility constraints, meaning that they deactivate the arcs between incompatible customers and vehicles (drivers). Constraints \eqref{time} help determine the service start time at every node on every vehicle while \eqref{rechargingtime} accounts for the time spent in a recharging visit. Constraints \eqref{timewindows} ensure that time windows are respected for every node. Additionally, \eqref{time}-\eqref{timewindows} also serve as sub-tour elimination constraints. Constraints \eqref{cap} guarantee that the demand of all customers is satisfied without exceeding vehicle capacity. Constraints \eqref{maxcap} ensure that the initial cargo level is non-negative and below the vehicle capacity. 
Constraints \eqref{soc} compute the \ac{SOC} at all customer nodes and ensure that it is not negative and does not exceed the battery capacity. Constraint sets \eqref{rechargingsoc} and \eqref{maxsoc} guide the calculation of the \ac{SOC} post charging, allowing partial recharging while \eqref{initial_icev} sets the energy range of \acp{ICEV}. 
Finally, the domains of decision variables are defined in  \eqref{binary} and \eqref{nonnegative}. In this formulation, the set $K$ may include \acp{EV} and \acp{ICEV} with different characteristics of which the subset $K^E$ is the set of \acp{EV}. 

\small
\begin{equation}
\mathcal{C}(\bm{\omega}, \bm{\lambda},T) =    \min  \left [ \sum\limits_{k \in K} f_k z_k + \sum\limits_{(i,j) \in A} \sum\limits_{k \in K} \big ( c_k \ p_{ijk}(q_{jk}, T) + M_k v_k  \big) \ t_{ijk} \ x_{ijk} + \sum\limits_{i \in N} g_i u_i  \right], 
\label{objective}
\end{equation}
\begin{align}
  \text{subject to: } 
  \sum\limits_{j:(i,j) \in A} \sum\limits_{k \in K} x_{ijk} &= 1-u_i & \forall i \in N \label{getunserved} \\
    \sum\limits_{i \in F} \sum\limits_{j:(i,j) \in A}  x_{ijk} &\leq 1 & \forall k \in K^{E}\label{atmost1} \\
    \sum\limits_{j:(i,j) \in A} x_{ijk} - \sum\limits_{j:(j,i) \in A} x_{jik} &= 0 & \forall i \in N', k \in K^{E} \label{flowEV} \\
    \sum\limits_{j:(i,j) \in A, \ j\not\in F} x_{ijk} - \sum\limits_{j:(j,i) \in A, \ j\not\in F} x_{jik} &= 0 & \forall i \in N, k \in K \setminus K^{E} \label{flowICEV} \\
    x_{ijk} + x_{jik} &= 0 & \forall i \in N_0, j \in F, k \in K \setminus K^E \label{NoFvisitICEV}\\ 
    \sum\limits_{j:(0,j) \in A} x_{0,j,k} = \sum\limits_{i:(i,n+1) \in A} x_{i,n+1,k} &= z_k & \forall k \in K \label{getvehicles} \\
    \sum\limits_{j:(i,j) \in A} x_{ijk} &\leq a_{ik} & \forall i \in N, k \in K \label{comp} \\
    \tau_{ik} + (s_i + t_{ijk}) - l_{n+1} (1 - x_{ijk}) &\leq \tau_{jk} & \forall i \in N \cup \{0\}, j \in N'_{n+1}, k \in K \label{time} \\
    \tau_{ik} + t_{ijk} + b_i (Y_{ik} - y_{ik}) - l_{n+1} (1 - x_{ijk}) &\leq \tau_{jk} & \forall i \in F, j \in N'_{n+1}, k \in K^{E} \label{rechargingtime} \\
    e_i \leq \tau_{ik} &\leq l_i & \forall i \in N'_{0,n+1}, k \in K \label{timewindows} \\
    q_{ik} - d_i + Q_k (1 - x_{ijk}) &\geq q_{jk} & \forall i \in N'_0, j \in N'_{n+1}, k \in K \label{cap} \\
    Q_k &\geq q_{0k} & \forall k \in K \label{maxcap} \\
    y_{ik}-p_{ijk}(q_{jk}, T)\ t_{ijk} x_{ijk} + E_k (1 - x_{ijk}) &\geq y_{jk} & \forall i \in N, j \in N'_{n+1} \label{soc}, k \in K \\
    Y_{ik} - p_{ijk}(q_{jk}, T) t_{ijk} x_{ijk} + E_k (1 - x_{ijk}) &\geq y_{jk} & \forall i \in F\cup\{0\}, j \in N'_{n+1}, k \in K^{E} \label{rechargingsoc} \\
    E_k &\geq Y_{ik} & \forall i \in F\cup \{0\}, k \in K^{E} \label{maxsoc} \\
    E_k &\geq y_{0k} & \forall k \in K \setminus K^{E} \label{initial_icev} \\
    x_{ijk}, z_k, u_i &\in \{0, 1\} & \forall i, j \in N'_{0,n+1}, k \in K \label{binary} \\
    y_{ik}, \tau_{ik}, q_{ik}, Y_{ik} &\geq 0 & \forall i \in N'_{0,n+1}, k \in K \label{nonnegative}
\end{align}
\normalsize

We note that each operational problem is a vehicle routing problem with compatibility constraints (which is NP hard \cite{Pillac2013AProblem}) along with capacity and electric range constraints. It becomes intractable to solve the operational problem exactly for realistic customer set sizes. In \Cref{methodology}, we resort to a state-of-the-art heuristic method, namely \ac{ALNS}, to solve the operational vehicle routing problem to determine the operational decisions on each day, given a strategic decision $K$. An exact solution of the strategic problem requires solving the operational problem for  all possible combinations of fleet mix, ambient temperature, and customer request set. If there is a cap on the total fleet size and/or number of vehicles of each kind in the mix, the number of fleet mix configurations would be finite, albeit numerous. However, it is impossible to enumerate all possible values of $T$ as it belongs to a continuous domain. Even if the temperature domain is discretized, enumerating all possible combinations of the customer set still remains a computationally intensive task. 
Thus, it is computationally intractable to determine an exact solution of the strategic problem. In \Cref{methodology}, we propose a novel heuristic solution approach based on \acf{SAA} \cite{Ahmed2002TheRecourse} to circumvent the tedium of exact computation.

\subsection{Power consumption model with temperature dependence}
 We propose an upgraded power consumption model enabling temperature dependence. We account for the following components of power consumption:
\begin{itemize} \setlength\itemsep{0em} \vspace{-2ex}
    \item load dependent mechanical power $P_M$ (following \citet{Bektas2011TheProblem},  \citet{Demir2011ATransportation}, and \citet{Goeke2015RoutingVehicles}),
    \item auxiliary power for temperature control $P_T$ (heating and air conditioning), and
    \item auxiliary power usage $P_O$. 
\end{itemize} 

\begin{figure}[H]
\begin{center}
\includegraphics[width=0.9\textwidth]{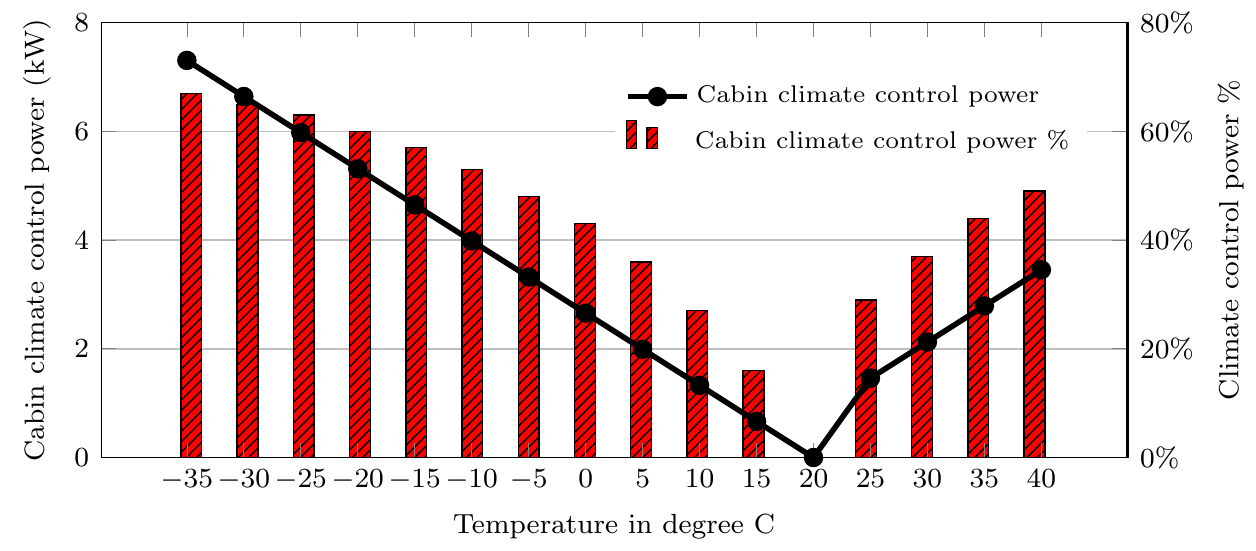}
\caption{Cabin climate control power variation for an electric mid-sized pickup van determined using the upgraded model for a desired cabin temperature of 20 \degree C: The left side axis tracks the power consumed for cabin climate control while the right side axis captures the climate control power as a percentage of the total power consumed by the vehicle}
\label{temp_power_kangoo}
\end{center}
\end{figure}

Ambient temperature has a significant effect on the cabin climate control energy consumption, and extreme temperatures can reduce greatly the range of electric cars. The extent of the impact can be seen through \Cref{temp_power_kangoo}, which demonstrates that for a mid-sized van, cabin climate control power is a significant component of the total vehicle power; at -10\degree C, it is about 4 kW, drawing about 53\% of the total vehicle power. The resultant overestimation of driving range by omitting cabin climate control power (auxiliary power) can threaten the feasibility of the implementation of electric mobility in urban logistics operations, and therefore it is important to capture auxiliary power in energy calculations.
Our model for the computation of auxiliary power uses the heat balance method as seen in \citet{Fayazbakhsh2013ComprehensiveMethod} and \citet{Valentina2014HVACVehicles}. The third component of auxiliary power usage refers to the power drawn from the car's energy source for running an external process such as maintaining a climate control box at a specific temperature in the case of the collection of blood samples from clinics (it excludes the power drawn by music systems and wipers since they depend on a separate 12V battery usually). Including the latter two new components of power has an effect equivalent to boosting all the cost coefficients of $x_{ijk}$'s in \eqref{objective}.

\tikzstyle{rect} = [draw, rectangle, fill=white!20, text centered, minimum height = 2em]
\tikzstyle{elli} = [draw, rectangle, fill=white!20, minimum height = 2em]
\tikzstyle{circ} = [draw, circle, fill=white!20,minimum width = 8pt, text centered, inner sep = 10pt]
\tikzstyle{diam} = [draw, diamond, fill=white!20,text width = 6em, text badly centered, inner sep = 0pt]
\tikzstyle{line} = [draw,-latex]

\begin{figure}[ht!]
\footnotesize

\begin{subfigure}{0.4 \textwidth}
    \begin{center}
    \begin{tikzpicture}[shorten >=1pt, node distance = 2cm, auto]
    \node [rect, fill = orange!40, rounded corners, text width = 7em, node distance = 3cm] (step1a) {Mechanical power};
    \node [rect, fill = yellow!50, rounded corners, text width = 8em, node distance = 2cm, below of = step1a] (step2) {Vehicle power};
    \path [line] (step1a) -- (step2);
    \end{tikzpicture}
    \caption{\citet{Goeke2015RoutingVehicles}}
    \end{center}
    \end{subfigure}%
    \begin{subfigure}{0.6 \textwidth}
    \begin{center}
    \begin{tikzpicture}[shorten >=1pt, node distance = 2cm, auto]
    \node [rect, fill = orange!40, rounded corners, text width = 7em, node distance = 3cm] (step1a) {Mechanical power};
    \node [rect, fill = blue!30, rounded corners, text width = 7em, node distance = 3cm, right of = step1a] (step1b) {Cabin climate control power};
    \node [rect, fill = green!30, rounded corners, text width = 7em, node distance = 3cm, right of = step1b] (step1c) {Auxiliary power};
    \node [rect, fill = yellow!50, rounded corners, text width = 8em, node distance = 2cm, below of = step1b] (step2) {Vehicle power};
    \path [line] (step1a) -- (step2);
    \path [line] (step1b) -- (step2);
    \path [line] (step1c) -- (step2);
    \end{tikzpicture}
    \caption{Current work}
    \end{center}
    \end{subfigure}
    \caption{Upgrade to the power consumption model}
    \label{f1}
\end{figure}
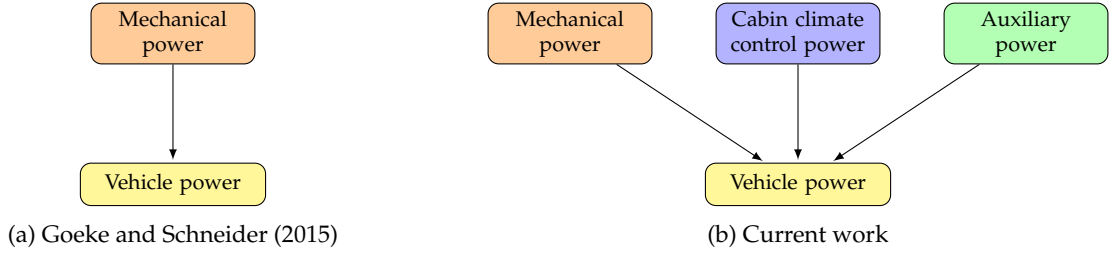
\normalsize

Power consumption $p_{ijk}(q_j, T)$ is a function dependent on the starting load at node $j$ and the ambient temperature $T$, when traveling from node $i$ to node $j$ by vehicle $k$. The following expression results in power in kW when using the parameter values indicated in Tables \ref{param_power_veh_type}, \ref{param_power_Pm} and \ref{param_power_Pt} in the appendix. 
\begin{align*}
p_{ijk}(q_j, T) &= \begin{cases}\phi^d \varphi^d P_{ijk}(q_j, T) & P_{ijk}(q_j, T) \geq 0, k \in K^E
\\ \phi^r \varphi^r P_{ijk}(q_j, T)  & P_{ijk}(q_j , T) <0, k \in K^E
\\  \max \left (kN'D + \frac{P_{ijk}(q_j)}{\eta' \cdot \eta_{tf}}, 0\right)  &  k \not \in K^E
\end{cases}
\\ P_{ijk}(q_j, T) &= P_M((i,j),k) +P_{T} ((i,j),k,T) + P_O 
\\ P_M((i,j),k) &=  0.001\times\left [c_r m_k(q_j) g \cos \alpha +\frac{1}{2} \rho_a A_f c_d v^2 +  m_k(q_j)g \sin \alpha + m_k(q_j)a \right ]v
\\ P_{T} ((i,j),k,T) & = \mathcal{I}_H \frac{1}{\eta_H} (P_H + P_{Con}) + \mathcal{I}_C \frac{1}{\eta_C} ( P_C  + P_{Con} + P_{Rad} + P_P),
\end{align*}
where $ P_M((i,j),k)$, $P_{T} ((i,j),k,T)$, and $P_O$ are mechanical, temperature control, and auxiliary power usage components respectively. 
The terms $\mathcal{I}_H$ and $\mathcal{I}_C$ are binary variables that indicate whether the heating mode or the cooling mode is being used in the definition of $P_T$. The components of $P_T$ are further defined below. 
\begin{eqnarray*}
    P_H &=& \dot{m} \ C_P \ (T_{d} - T) \\
    P_C &=& \dot{m} \ C_P \ (T-T_{d}) 
    \\ P_{Con} &=&  \sum_{i \in Surfaces} A_i U_i \abs{T  - T_{d}} \\
    P_{Rad} &=& \sum_{i\in Surfaces} A_i \tau_i {I}_{DN} \bigg ( \cos \theta_{i}+ C \frac{(1+ \cos\Sigma)}{2}  +\rho_g \cos \theta_z \frac{(1-\cos \Sigma)}{2} \bigg),\\
    \text{where } I_{DN} &=& A \exp (-B/\cos \theta_z) \\ 
    P_P &=& \sum_{i \in Persons} H_{Pr} A_{Du} 
\end{eqnarray*}
The main component is the ventilation load, i.e., the power required to warm up ($P_H$) or cool down ($P_C$) the air entering the car. The rest of components are the following: $P_{Con}$ is the power loss due to conduction through the surfaces of the vehicle and/or exterior and interior convection, $P_{Rad}$ represents the cooling load introduced in the car due to direct and diffuse solar radiation, and $P_P$ is the cooling load corresponding to the metabolic heat generated by humans seated in the vehicle.
A brief description of all the parameters along with pointers to their sources is presented in tables \ref{param_power_veh_type}, \ref{param_power_Pm} and \ref{param_power_Pt} in the appendix. Specifically, the efficiency parameters used in the definition of $p_{ijk}(q_j, T_o)$ are defined in \Cref{param_power_veh_type}, and the parameters used in the definitions of $P_M$ and $P_T$ are defined in \Cref{param_power_Pm} and \Cref{param_power_Pt}, respectively.

\section{Methodology\label{methodology}}
Since the \ac{SFSMP} is a two-stage stochastic program, we employ a fairly standard method, known as \acf{SAA}, used in the literature for pursuing solutions of two-stage stochastic programs. 
We propose solving the \ac{SFSMP} by \ac{SAA}  through Monte Carlo sampling of the requests and ambient temperature using known probability distributions. 

\subsection{Solving the strategic problem}
We propose a nested framework of solving the fleet sizing problem with lack of knowledge on the occurrence of demand requests and daily temperatures at the strategic level. The framework is
is delineated through a flow chart in \Cref{f1}. For every initial fleet mix, the set of customer requests and the ambient temperature are sampled. On the resultant operational instance, the \ac{ALNS} algorithm is applied to  determine the operational cost for each sample instance $i$, which is $\mathcal{C}(\bm{\omega},, \bm{\lambda}^i,T^i) $. The half width of the confidence interval for the sample average estimate based on the samples generated so far is determined. If it meets the required criterion (for e.g., of falling under a certain threshold), the process proceeds to computing the sample average operational cost for the considered fleet mix. In this step, we approximate the expectation in \eqref{strat_obj} with a sample average of costs of operation. 
\begin{equation}
\mathbb{E}_{\mathcal{N}, \mathcal{T}} \left [ \mathcal{C}(\bm{\omega}, \bm{r},\mathcal{T})  \right ] \approx \frac{1}{n} \sum_{i =1}^n \mathcal{C}(\bm{\omega}, \bm{\lambda}^i,T^i) , \label{strat_approx_obj} 
\end{equation}
where $n$ is the number of samples with the simulated outcomes of customer set determined by $\bm{\lambda}^i$ and ambient temperature $T^i$,  for each operational period $i$. After evaluating the sample average operational cost for all the fleet mixes of interest, the best fleet mix is determined. 

 We remark that an exhaustive enumeration of fleet mixes may result in an exponentially large number of fleet mixes if the master list of vehicles is heterogeneous and large. Thus, it is suggested that the set of fleet mixes may be determined heuristically or exhaustively, depending on the size of the problem.  Given the list of preferred commercial vehicles for inclusion in the fleet, enumeration with reasonable step sizes would be practical. 

\tikzstyle{rect} = [draw, rectangle, fill=white!20, text centered, minimum height = 2em]
\tikzstyle{elli} = [draw, rectangle, fill=white!20, minimum height = 2em]
\tikzstyle{circ} = [draw, circle, fill=white!20,minimum width = 8pt, text centered, inner sep = 10pt]
\tikzstyle{diam} = [draw, diamond, aspect = 6,  fill=white!20,text width = 6em, text badly centered, inner sep = 0pt]
\tikzstyle{line} = [draw,-latex]

\begin{figure}[ht!]
\footnotesize
\begin{center}
    \begin{tikzpicture}[shorten >=1pt, node distance = 1cm, auto]
    \node [rect, fill = blue!30, rounded corners, text width = 12em, node distance = 2.5cm] (step1) {Set initial fleet mix  using membership vector $\bm{\omega}$};
    \node [draw = none, below of = step1, node distance = 1cm] (step1ba) {};
    \node [rect, fill = blue!30, rounded corners, text width = 10em, node distance = 2.5cm, right = 3.8cm of step1ba] (step1b) {Generate next fleet mix vector $\bm{\omega}$};
    \node [rect, fill = orange!30, rounded corners,  below of = step1, text width = 13em, node distance = 3cm] (step2) {Sample customer request vector $\bm{\lambda}$ \\ and temperature $T \in \Omega (\mathcal{T})$};
    \node [draw = none, right = 2.1cm of step2] (step2b) {};
    \node [rect, fill = green!50, rounded corners, below of = step2, text width = 12em, node distance = 2cm] (step3) {Apply \ac{ALNS} ($\bm{\omega}$, $\bm{\lambda}$, $T$)};
     \node [diam, fill = orange!30, rounded corners, below of = step3, text width = 20em, node distance = 2cm] (step4) {Is half width criterion satisfied?};
     \node [rect, fill = yellow!50, rounded corners, below of = step4, text width = 12em, node distance = 2cm] (step4a) {Compute sample average operational cost ($\bm{\omega}$)};
     \node [diam, fill = blue!30, rounded corners, below of = step4a, text width = 20em, node distance = 2cm] (step5) {Are all the mixes of interest evaluated?};
     \node [rect, fill = red!30, rounded corners, below of = step5, text width = 12em, node distance = 2.3cm] (step6) {Determine the best fleet mix};
    \path [line] (step1) -- (step2);
     \path [line] (step2) -- (step3);
     \path [line] (step3) -- (step4);
     \path [line] (step4) -- (step4a);
     \path [line] (step4a) -- (step5);
     \draw[draw = none] (step4) -- node[right, near start, pos = 0.3] {Yes} (step4a);
     \path [line] (step2b.center) -- (step2);
     \draw (step4) -| node[above, near start, pos = 0.1] {No} (step2b.center);
     \path [line] (step1b.west) -- (step1ba);
     \draw (step5) -| node[above, near start, pos = 0.1] {No} (step1b.south);
     \path [line] (step5) -- (step6);
     \draw[draw = none] (step5) -> node[right, near start] {Yes} (step6);
     \draw[fill = yellow!10, opacity = 0.1,dashed, color = black!40] (-4.1,-1.8) rectangle (4.9,-7.8);
     \node[draw=none] at (-2.3,-1.5) {\textcolor{gray}{Simulation of operational days}};
    \end{tikzpicture}
    \caption{Solving the \acf{SFSMP} using \acf{SAA}}
    \label{flow_SAA}
\end{center}
\end{figure}
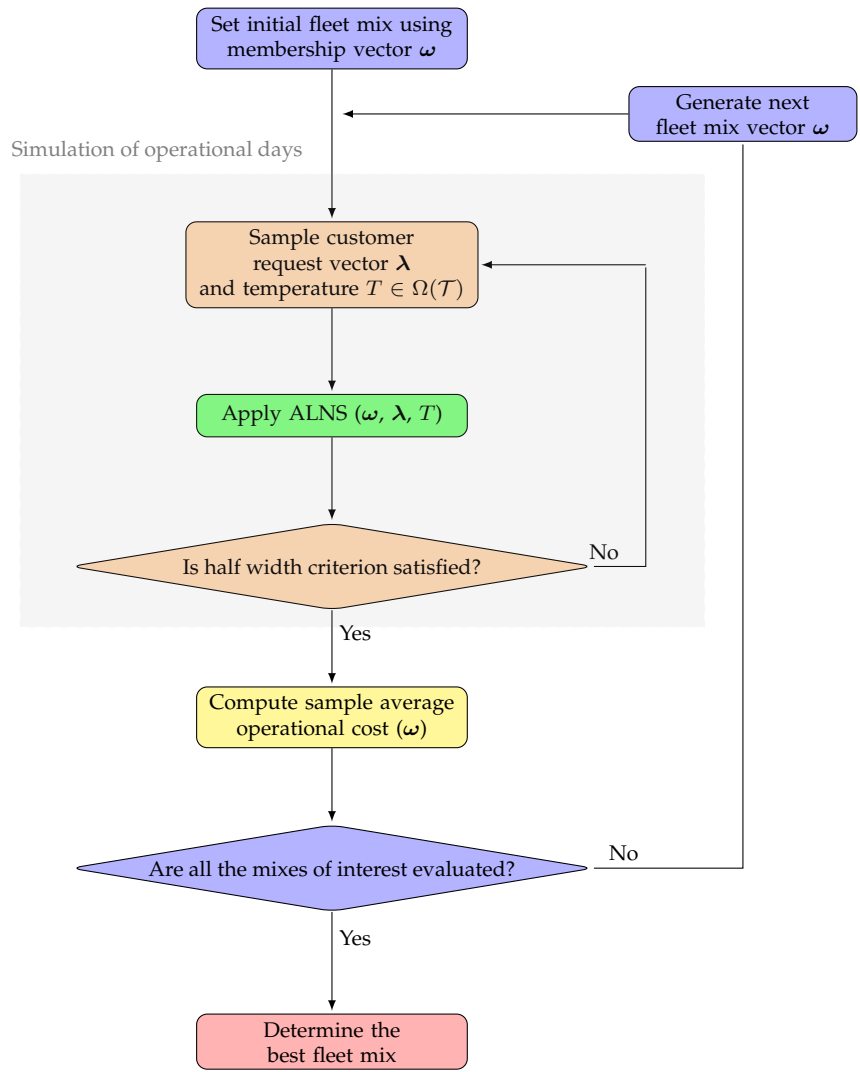

\normalsize

\subsection{Solving the operational problem}
In this section, we present the solution method used to determine decisions that must be implemented in each period of operation. 
Since determining an exact solution of the operational problem presented in \Cref{prob_form} is not computationally tractable for realistically sized instances, we pursue a heuristic approach. In particular, we develop a \ac{ALNS} algorithm (\citet{Rpke2006AnWindows}), which has been utilized extensively in the literature to solve vehicle routing problems and specifically, the mixed fleet vehicle routing problem with time windows, as can be seen in \citet{Goeke2015RoutingVehicles} and \citet{Keskin2016PartialWindows}.
We use the same parameters as those used in \citep{Christensen2014AdaptiveThesis}.

 We have employed the following  destroy operators for \textit{customer removal} in the algorithm:
\begin{enumerate}
\item Random removal \citep{Rpke2006AnWindows}: A given number of  randomly chosen customers are removed from the tours.
\item Worst removal \citep{Rpke2006AnWindows}: A given number of customers that result in the highest insertion costs are removed.
\item Shaw removal \citep{Rpke2006AnWindows}: A given number of customers related to each other through a randomly selected initial customer according to five different criteria (customer proximity, tour proximity, demand difference, arrival time difference and time window similarity) are removed.
\item Random route removal: Entire tours are removed from a given solution until a given number of customers are removed. 
\item \ac{SISR} \citep{ChristiaensJanSIbS}: This destroy operator was introduced in  \citep{ChristiaensJanSIbS}. Strings of nodes or pairs of strings of nodes separated by a few nodes are removed. 
\end{enumerate}
The following repair methods have been used for \textit{inserting the removed customers}:
\begin{enumerate}
    \item Greedy insertion \citep{Rpke2006AnWindows}: The customer node to be inserted and the position of insertion are determined by minimizing the cost of insertion among all feasible insertions. 
    \item 2-Regret insertion \citep{Rpke2006AnWindows}: In this repair method, we select the customer for which the difference between the second lowest insertion cost and the lowest insertion cost is maximum. This customer is then inserted in its best position.    
\end{enumerate}
The insertion methods are modified to allow at most one recharging visit on each tour. In this \ac{ALNS}, a simulated annealing based acceptance criterion is employed to guide the search. The start temperature is calculated such that a solution $s$ times worse than the starting solution is accepted with a probability of 50\%, in the first iteration.
\begin{align*}
    T_{Start} = - \frac{sf(x)}{\text{ln}(0.5)}
\end{align*}
Additionally, a resetting mechanism is used to refresh the current solution with the best solution after a certain number ($\gamma$) of iterations without any improvement to the best solution. We present the other parameters used in the \ac{ALNS} in \Cref{param_ALNS} of \Cref{app_param}.

\section{Case Studies \label{cases}}
As part of the EU-funded research project called Electric Urban Freight and Logistics (EUFAL), the  aim of which is to encourage the electrification of commercial fleets in order to reduce the carbon footprint of urban logistics operations, we developed two real case studies. We present our analysis of the two cases  in this section. 

\subsection{Blood sample collection from private clinics for Region Hovedstaden}
Region Hovedstaden (Region H: \url{www.regionh.dk/}) is a public authority that manages healthcare and social services in the Capital Region of Denmark, which corresponds roughly to the metropolitan area of Copenhagen.  Collecting blood samples from doctors and delivering these to a testing facility is one of the functions carried out by the Region, which is the focal logistics process of this case study. Everyday, blood samples must be collected from the clinics of private physicians and delivered to the laboratory testing facility in a hospital by vehicles that start their trips from a different hospital that is about four kilometers away from the first hospital. Due to the perishability of the collected samples, a 150 minute shift is designed for all drivers. Since the viability of the samples degrades beyond three hours after extraction, physicians would be able to send the samples by taxi if the fleet operation plan does not pick up from the corresponding physician location. On each work day, the Region operates two identical blood sample collection shifts - one in the forenoon and one in the afternoon. The Region aims to optimize its mix of delivery fleet and make smart operational decisions to minimize the cost of missed doctor visits and the cost of  energy used for routing. 
We consider the problem of strategic sizing of the fleet required to manage these operations. Similar operational problems were considered by \citet{Grasas2014} in Catalonia, \citet{Zufferey2016} in Geneva, and \citet{Anaya-Arenas2014} and  \citet{Kergosien2014AServices} in Quebec. 

In line with its mission, the Region plans to use an all-electric fleet for its logistics operations. It has access to two different types of vehicles, medium-sized electric pickup vans and 
electric cargo bikes 
 with driving ranges of 270 and 140 kilometers respectively.  The characteristic parameters of the two types of vehicles considered in this study are presented in \Cref{blood_veh_features} and are inspired from those of the mid-sized electric pickup van Renault Kangoo and the electric cargo bike TRIPL. Though the cargo bike does not have to manage cabin climate control, the climate box used for maintaining the temperature of the collected blood samples at 20 degree Celsius draws energy from its battery directly with a base power of 0.1 kW in addition to the actual power required to heat or cool the box. In the electric van, the climate box is hooked to its secondary 12V battery, thus not affecting the driving range of the vehicle.    
 
 \begin{table}[ht!]
    \centering
\caption{Region H case: Features of the vehicles considered}
\label{blood_veh_features}
    \begin{tabular}{rll} 
    \toprule
    \textbf{Feature} & \textbf{Pickup Van} & \textbf{Cargo bike} \\
    \midrule
    Kerb weight &	1426 kg & 322 kg \\  
    Max speed &	130 kmph &	45 kmph \\ 
    Operational speed &	45 kmph &	40 kmph \\ 
    Battery &	33 kWh &	8.64 kWh \\ 
    Capacity &	720 vials &	288 vials \\
     Additional weight  &	200 kg &	50 kg \\
    Acquisition cost &	28,223 USD  &	17,900 USD \\
    Auxiliary power & 0 & 0.1 kW
    \\
     \bottomrule
    \end{tabular}
\end{table}
 
\begin{figure}[h!]
\begin{center}
\includegraphics[scale = 0.9]{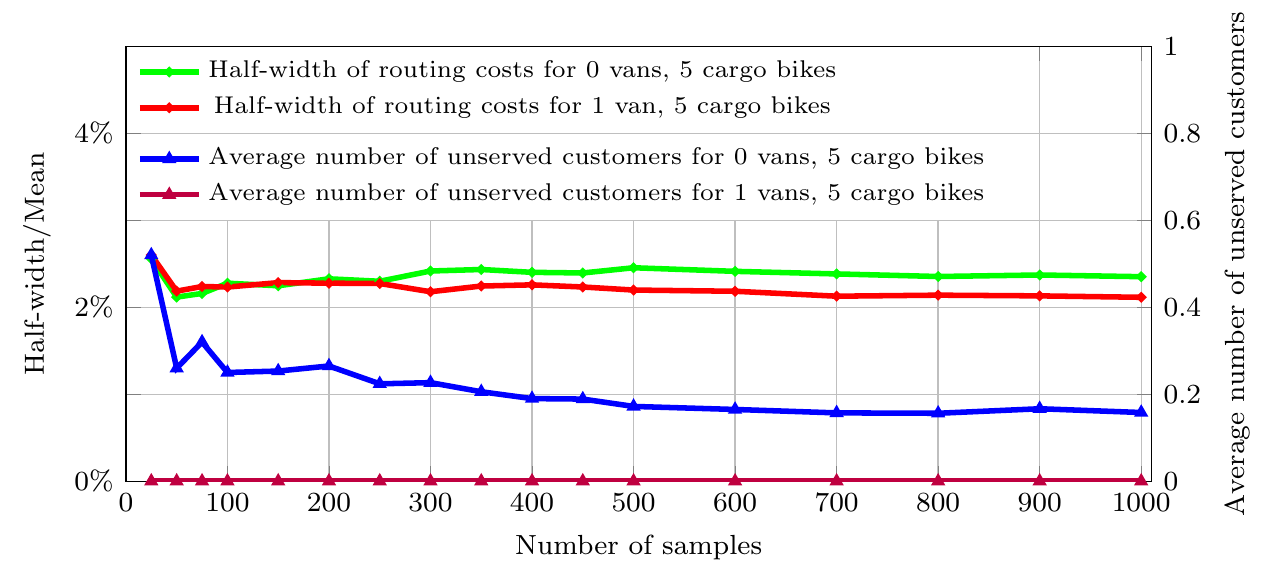}
\caption{Region H case: Half width analysis for the fleet mix containing [0 vans, 5 cargo bikes] and [1 van, 5 cargo bikes]}
\label{blood_hw}
\end{center}
\end{figure}

The Region is now faced with the decision of determining the number of vehicles of each type required to constitute its fleet for its long term operations. At the operational level, the requests from doctors are revealed on the day of operation. However, the actual demands are not known days in advance.  Based on data collected through a survey we conducted, we found that the demands are seasonal and can be assumed to follow different distributions in \textit{summer} and \textit{non-summer} parts of the year. A detailed description of the generation of operational instances for this case study may be seen in \Cref{blood_inst_gen}.  Based on initial experiments, we have found that the maximum fleet size will never exceed $15$. Thus, we consider all  vehicle mix possibilities such that the total fleet size does not exceed $15$. A horizon of 10.6 years  with 227 working days  per year and two work shifts per day is considered as the horizon based on the average lifetime of vehicles stated in \citet{Mitropoulos2017TotalVehicle}.   For every possible fleet mix, an estimate of the operational cost is determined over the entire horizon. To this estimate, the acquisition cost of the corresponding fleet mix is added to determine the \ac{TCO}. The best fleet mix is thus determined using the procedure described in \Cref{flow_SAA}. Since we adopt \ac{SAA} to solve the \ac{SFSMP}, the approximation of the expected cost relies on a sample average, for computing which, choosing the right sample size is crucial. \Cref{blood_hw} presents an analysis of the variation of the half width of the 95\% confidence interval of the average routing cost for two fleet mixes at various sample sizes. We note that the half width stabilizes beyond 500 samples. However, considering the average number of unserved customers (or equivalently, the average cost of missing customers), we note that 600 or more samples would be required for the [0 vans, 5 cargo bikes] fleet mix. However, there are no unserved customers for the other mix. Considering both these aspects, 600 samples for the former and 500 samples for the latter mix would be sufficient to obtain a meaningful sample average. It is clear that the required sample size depends on the actual configuration considered and that the solutions of interest have a sample size requirement under 1000 samples. To obtain the following results, we use 1000 samples for all the fleet mixes. 

\begin{figure}[ht!]
\begin{center}
\includegraphics[scale = 0.9]{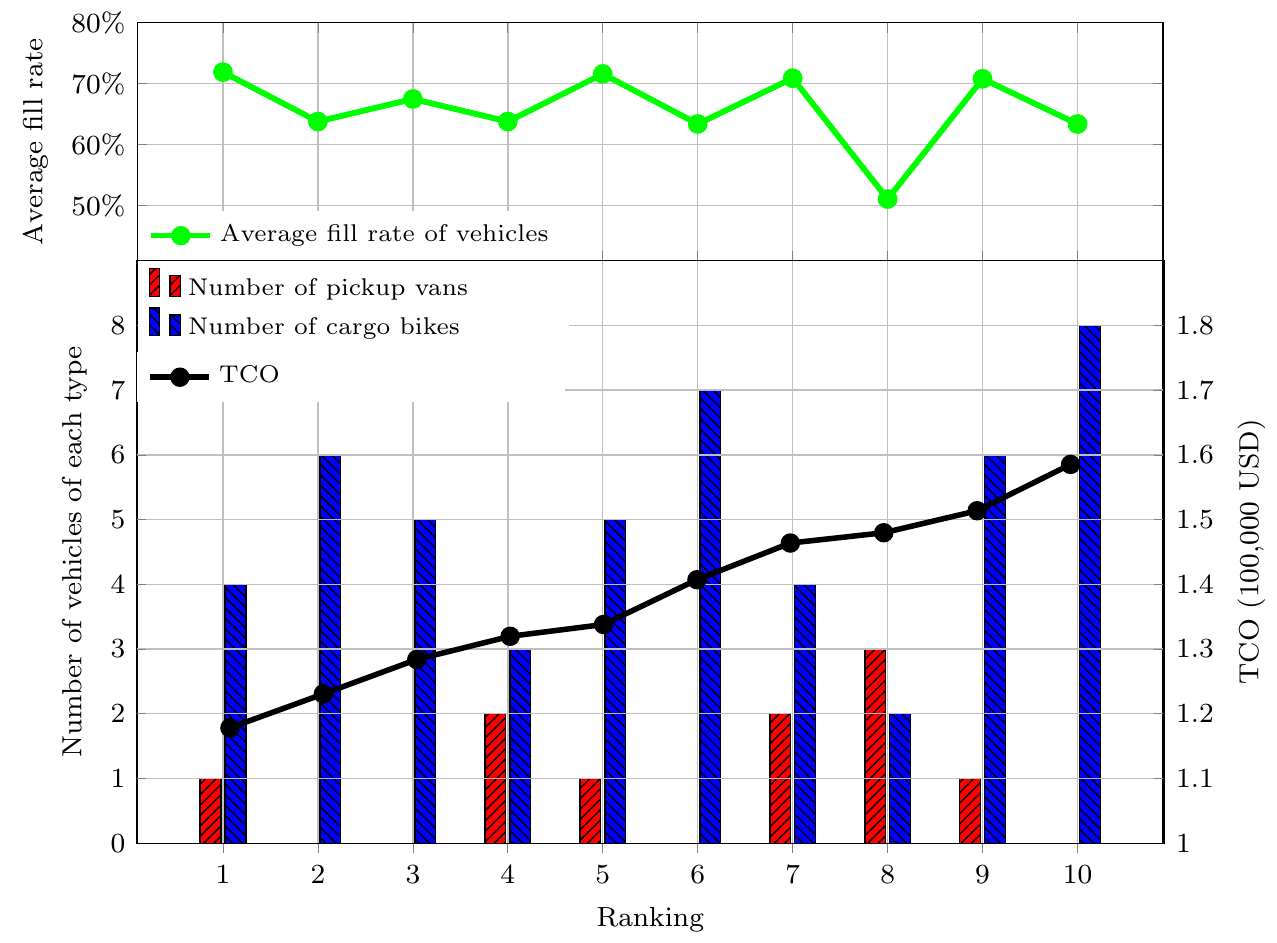} 
\caption{Region H case: The ten best fleet mix solutions determined by \ac{SAA} with their composition, \ac{TCO}, and average fill rates} \label{blood_top_ten_TCO}
\end{center}
\end{figure}

The fleet mix, the \ac{TCO}, and the average fill rate of the ten best heuristic solutions that minimize the \ac{TCO} of the fleet  are presented in \Cref{blood_top_ten_TCO}. We note that the tenth best solution to the stochastic problem is about $40,000$ US dollars costlier than the first best solution. Such a difference (25\% of the first best solution's \ac{TCO}) reaffirms the need to incorporate the impact of uncertainty into the study. The composition of the fleet mixes of the top ten solutions indicate that there is a strong preference for the inclusion of cargo bikes in the mix and that pickup vans, though not as popular, do appear in six positions. It may be noted that broadly, a pickup van  is equivalent to (and slightly more expensive than) two cargo bikes in the marginal sense. We note that in spite of the variation of the fleet mixes, the average capacity fill rates of nine of the ten solutions fall in the 50-80\% range. 

\begin{figure}[ht!]
\begin{center}
\includegraphics[scale = 0.8]{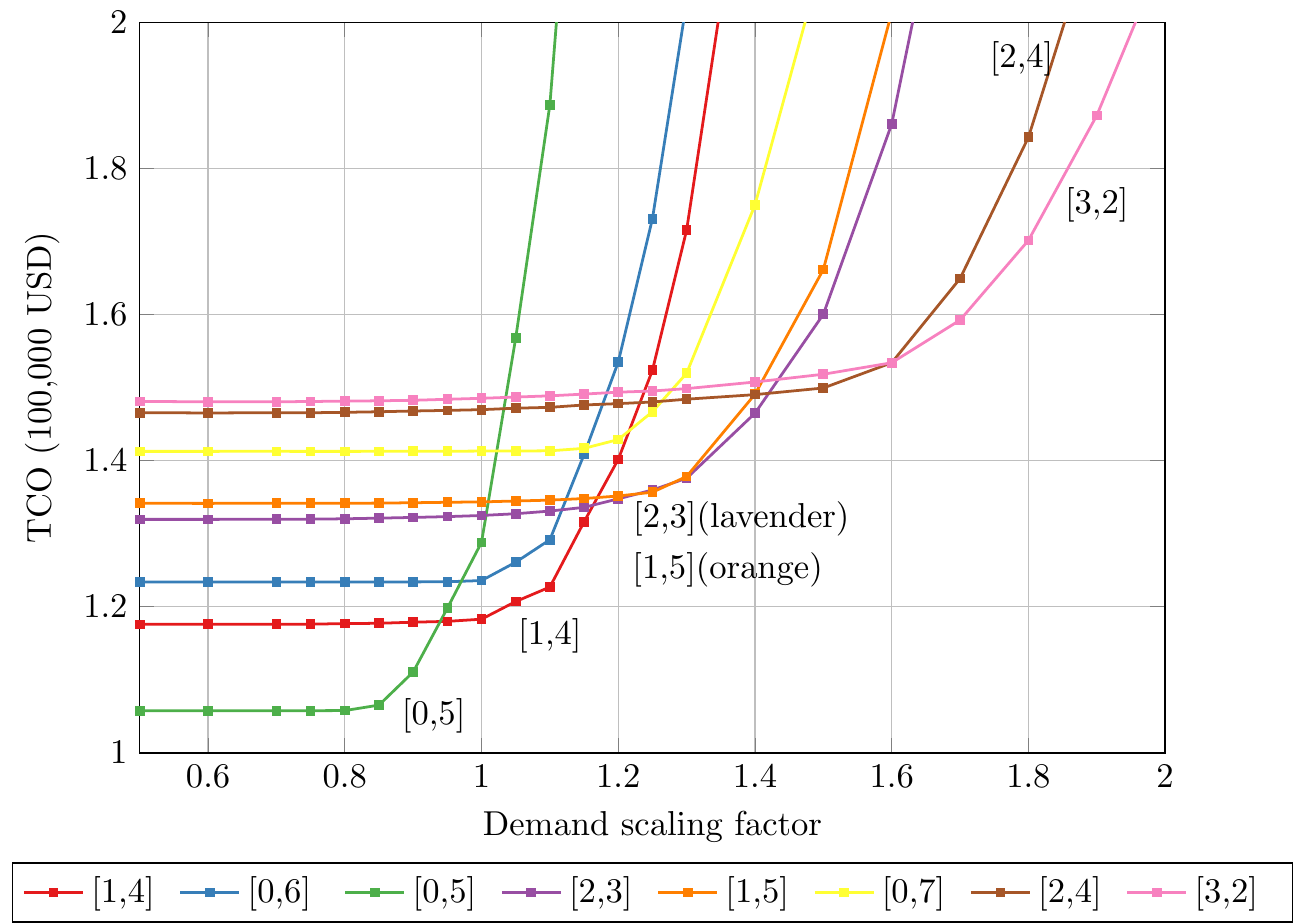}
\caption{Region H case: Sensitivity analysis of the top eight fleet mix solutions to variation of scaling demand at all clinics. In the legend, an entry with $[a,b]$ refers to the fleet mix with $a$ vans and $b$ cargo bikes.} \label{blood_demScaling}
\end{center}
\end{figure}

The best fleet mixes optimize for a given demand characterization, leaving little room for planned redundancy. While operating logistics systems, planners prefer adding a safety factor to build in resilience against near-future changes due to the overall evolution of the system. In \Cref{blood_demScaling}, we present a plot of the variation of \acp{TCO} with increasing level of demand scaling. The base case of the scaling factor set to 1 is the current demand situation faced by the Region. The remaining factors scale down or up the demand faced at every clinic. We note that the best mix [1 van, 4 cargo bikes] remains a preferred solution for the  interval [0.95, 1.15] of the scaling factor. However, when planning for scenarios with scaling between [1.2, 1.3], the mix [1 van, 5 cargo bikes] is preferred. We find that at higher \acp{TCO}, there is greater safety against fluctuations in demand levels. Although these levels of variation in demand may not seem very realistic for the current case study, this analysis would be applicable to other case studies dealing with different products.   Thus, a firm planning to add in redundancy into its fleet management may seek a fleet mix with higher \ac{TCO} to ensure readiness for overall higher demand scenarios. 

\begin{figure}[ht!]
\begin{center}
\includegraphics[scale = 0.9]{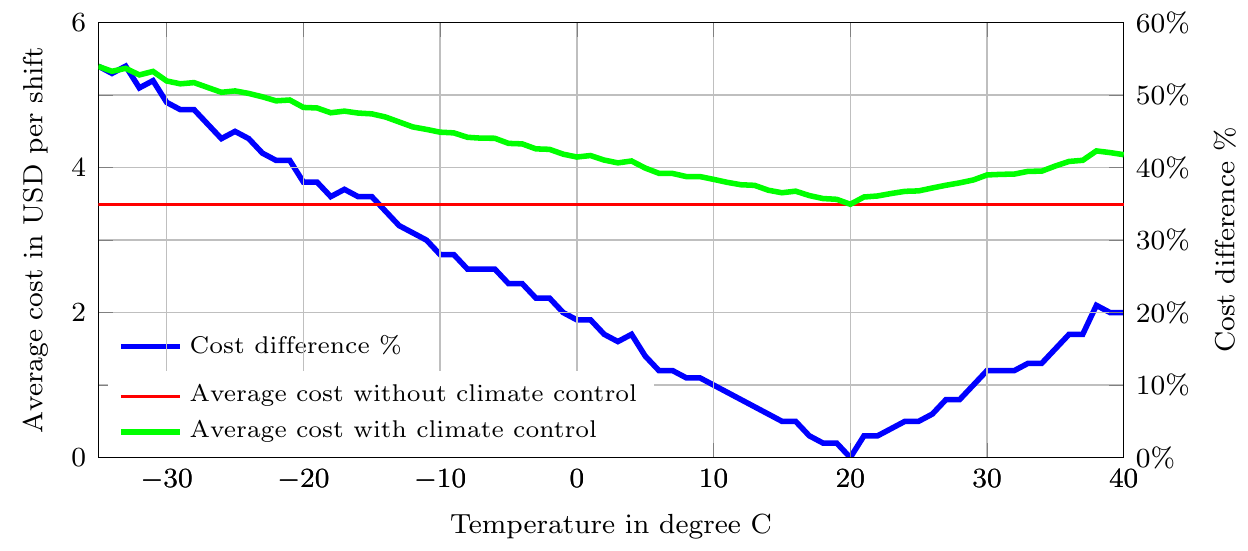}
\caption{Region H case: Average operational cost per shift in USD (left axis) and percentage cost difference over the case of no cabin climate control (right axis) with varying ambient temperature for the best fleet mix solution with 1 van and 4 cargo bikes}
\label{temp_1K4T}
\end{center}
\end{figure}

\Cref{temp_1K4T} presents the variation of average operational cost for the best fleet mix [1,4] in response to a variation of temperature. We note that sub-zero temperatures lead to a higher overall energy consumption compared to tropical temperatures. At $0$ degrees Celsius, not including the effect of temperature on vehicular power consumption leads to a 20\% difference from the base case (of 20 degree Celsius since the desired cabin temperature is 20 degrees Celsius) where the effect of temperature is not considered. This affects the service plan critically resulting in the mismanagement of the logistics system, since the available energy reserve as well as the driving range are overestimated during planning if the effect of cabin climate control and auxiliary energy usage are not considered.

Thus, the case study of the Region Hovedstaden reveals important insights regarding the importance of modeling uncertainty, the sensitivity of solutions to variations in demand levels, and the effect of temperature on planning fleet operations. 

\subsection{Technician routing for MT Højgaard}

In this case, a Danish construction company,  MT Højgaard (MTH), is currently considering electrification of the commercial fleet of its subsidiary Lindpro, provided there is no negative impact on their business performance. Lindpro is an electric installations company that provides electrician services to various customer locations. Each technician  has a specific skill level and technicians with higher skill levels are capable of attending to tasks of skill levels lower than their skill level. However, in reality, the company allows technicians of a specific skill level serving at most one level below their skill level. The company currently operates small and large \acp{ICEV} in its fleet. In this case study, we consider the strategic decision of investment in the fleet (as if the entire fleet must be acquired now). Only the small \acp{ICEV} are currently being considered for replacement with \acp{EV} as the tasks and routes operated by the large \acp{ICEV} are not yet relevant for replacement with \acp{EV}. Each driver begins and ends his/her tour at their homes. 

\citet{Villegas2018TheVehicle} considered a technician routing problem with mixed electric fleet and proposed a decomposition based solution method. We endeavor to determine the fleet mix using \ac{SAA}. The \ac{ALNS} is modified to allow compulsory visits to a driver's home at the start and end of their tours if the corresponding vehicle is used. 

We  generated a master list of customers in the Greater Copenhagen region with potential tasks that are associated with service times, time windows, and skill requirements inspired from historical information.  Compatibility matrices were developed with the following rules based on discussions with MTH:
\begin{itemize}
\item There are four skill levels that enable matching customer requests and technicians, where Level 1 is the most basic (handyman) while Level 4 is specialty (classified technician).
\item Technicians of a particular skill level step down by at most one level. 
\item It is fairly common for technician of a certain level to perform tasks one level below. However, it is very rare for Level 2 technicians to perform Level 1 tasks since Level 1 technicians belong to a union that is different from that of technicians of other levels. 
\end{itemize}

A detailed description of the various aspects of sample instance generation is provided in \Cref{mth_inst_gen}.

\begin{figure}[ht!]
\begin{center}
\includegraphics[scale = 0.9]{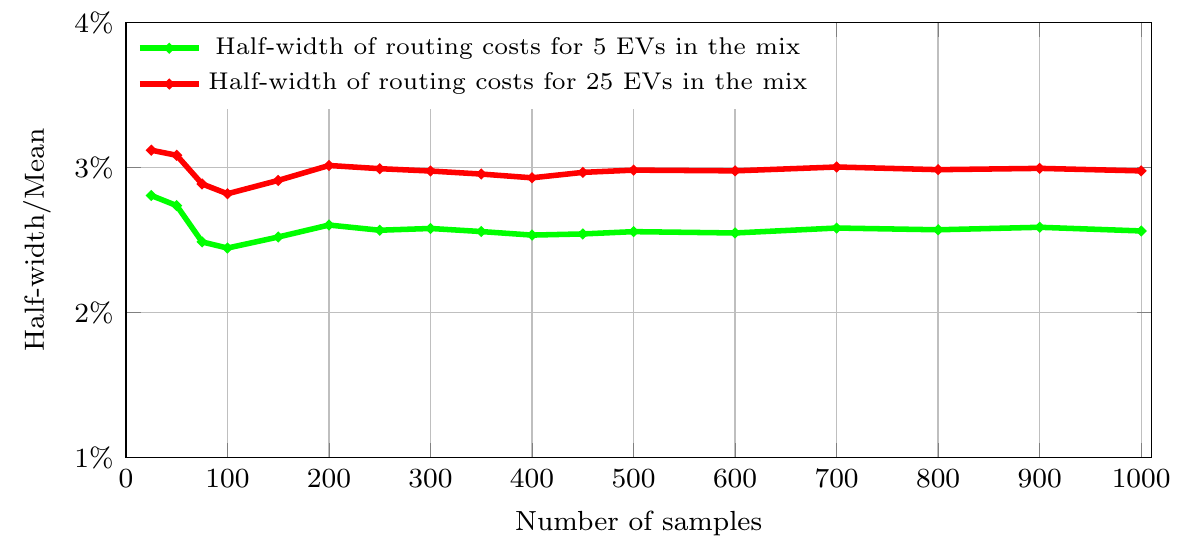}
\caption{MTH case: Half width analysis for the fleet mix containing 5 and 25 \acp{EV}}
\label{mth_hw}
\end{center}
\end{figure}

As part of implementing \ac{SAA}, we determined the suitable sample size based on the stabilization of the half width of the 95\% confidence interval of the sample average of the operational cost (see \Cref{mth_hw}). We found that  a sample size of 200 would be sufficient.

\begin{figure}[ht!]
\begin{center}
\includegraphics[width=0.8\textwidth]{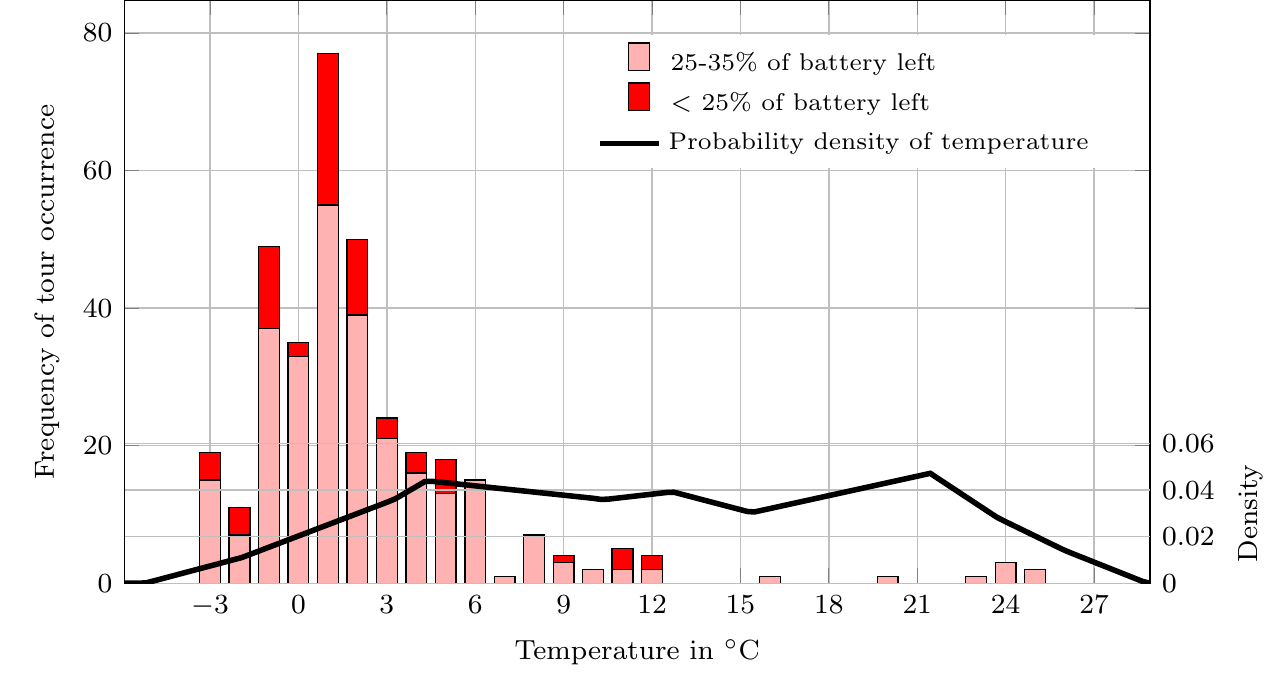}
\caption{MTH case: A stacked histogram (left axis) presents the frequency of occurrence of \ac{EV} tours in which the \ac{SOC} at the end of the tour falls between 25-35\% (pink) and below 25\% (red) based on computational results from 200 sample days with 110,786 \ac{EV} tours. The probability density function (a relative measure obtained by summing over the probability densities of four seasons) of temperature is  presented in black  (right axis).}
\label{mth_range_anxiety}
\end{center}
\end{figure}

\begin{figure}[ht!]
\begin{center}
\includegraphics[scale = 0.9]{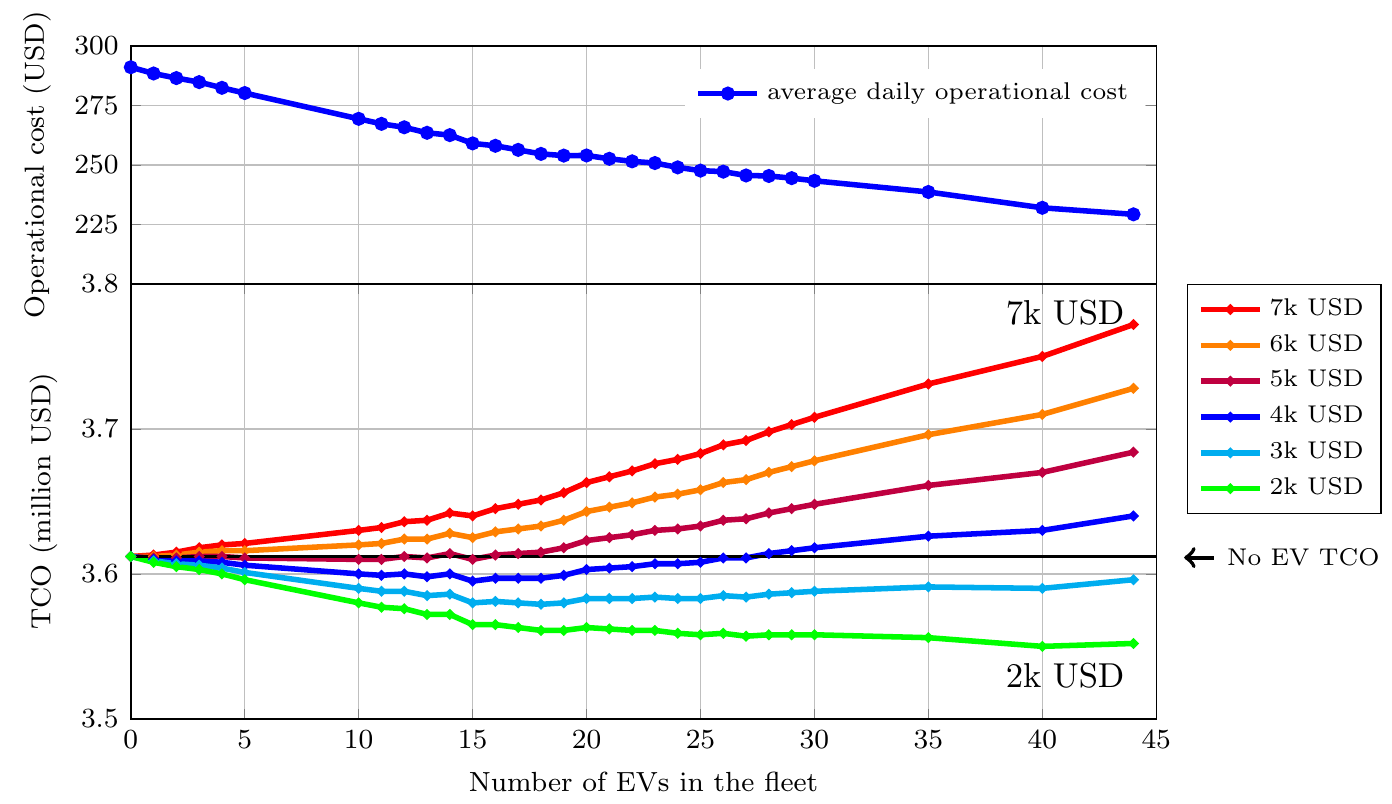}
\caption{MTH case: (a) Variation of average daily operational cost and (b) Variation of \ac{TCO} with \ac{EV} composition at various levels of price difference between \ac{EV} and  small \ac{ICEV} when the small \ac{ICEV}'s price is fixed at 21k USD, the price of diesel is 1.34 USD/liter, and the price of electricity is 0.1973 USD/kWh. } 
\label{mth_group}
\end{center}
\end{figure}

 \begin{table}[ht!]
    \centering
\caption{MTH case: Features of the vehicles considered}
\label{mth_veh_features}
    \begin{tabular}{rlll} 
    \toprule
    \textbf{Feature} & \textbf{\ac{EV}} & \textbf{Small \ac{ICEV} } & \textbf{\ac{ICEV}}\\
    \midrule
    Kerb weight &	1430 kg & 1623 kg & 1800 kg \\  
    Operational speed &	43.2 kmph &	43.2 kmph & 43.2 kmph \\ 
    Battery &	33 kWh & -	 & - \\ 
    Additional weight  &	250 kg &	250 kg & 400 kg \\ 
    Acquisition cost (USD) &	32,000 USD  &  21,000 USD & 28,000 USD	 \\	
    Auxiliary power & 0 & 0 & 0 
    \\     \bottomrule
    \end{tabular}
\end{table}
The features including the additional load carried as electrician equipment in the vehicles are presented in \Cref{mth_veh_features}. Since the 44 small vehicles can be either \acp{ICEV} or \acp{EV}, we considered fleet mixes with 0, 5, $\dots$, 44 \acp{EV} and with the remaining small vehicles being \acp{ICEV}. When 5 \acp{EV} are included in the fleet, the first 5 drivers out of the 44 drivers are assigned \acp{EV}. The list of drivers is in the order of the drivers' \ac{EV} eligibility, which is determined by the proximity of their home addresses to the nearest charging stations (see \Cref{mth_inst_gen}). On solving the \ac{ALNS} for each of these fleet mixes for each of the sample days, we determine the average daily operational cost  (see the upper plot of \Cref{mth_group}), which reduces with increasing number of \acp{EV} in the fleet. The \ac{ALNS} also generates the tours that must be executed on each day. \Cref{mth_tour_stats} presents some statistics pertaining to prescribed tours by vehicle type. On average the \acp{EV} in the fleet spend greater duration driving compared to the \acp{ICEV}. Correspondingly, the average task time (also considering breaks as tasks) of \acp{EV} is also more than an hour lower than that of \acp{ICEV} operated in the fleet. Thus, \acp{EV} are deployed to perform less time-consuming tasks and drive between customer locations more often since the energy consumed in driving \acp{EV} is cheaper than that for \acp{ICEV}. In spite of the longer duration spent on driving around, the average amount of energy consumed by \acp{EV} is 5.6 kWh based on results from 200 sample days. Though the planned energy usage of \acp{EV} very rarely exceeds 16 kWh on any tour, there are two occurrences when the \ac{SOC} at the end of the tour drops below 6\% of the usable battery capacity of 30 kWh, 24 occurrences when it drops below 20\%, and 148 occurrences when it drops below 30\%. These occurrences usually take place under 1 \degree C. \Cref{mth_range_anxiety} presents a measure of range anxiety, i.e.,  the temperature-dependent risk of the \ac{SOC} falling below 25\% and between 25-35\% of the battery capacity in terms of the number of times it happens in 200 sample days with 110,786 prescribed tours. We note that at low temperature, drivers run a significantly high risk of facing range anxiety. This observation indicates that if the power consumption model had not accounted for temperature-dependent power, such a prescription of tours would result in en-route driving range chaos. 

\begin{table}[ht!]
\centering
\caption{MTH case: Tour statistics by vehicle type}
\label{mth_tour_stats}
\begin{tabular}{r|ccccc}
\toprule
\textbf{Vehicle} &\begin{tabular}[c]{@{}c@{}} \textbf{Customers} \\ \textbf{served} \end{tabular} & \begin{tabular}[c]{@{}c@{}} \textbf{Driving} \\ \textbf{time (h)} \end{tabular} & \begin{tabular}[c]{@{}c@{}} \textbf{Serving} \\ \textbf{time (h)} \end{tabular} & \begin{tabular}[c]{@{}c@{}} \textbf{Avg. task} \\ \textbf{time (h)} \end{tabular} & \begin{tabular}[c]{@{}c@{}} \textbf{Tour cost} \\  \textbf{(USD)} \end{tabular} \\
\midrule
\ac{EV} & 7.86 & 1.05 & 8.47 & 1.54 & 2.77 \\
Small \ac{ICEV} & 4.11 & 0.42 & 8.52 & 2.70 & 2.31 \\
\ac{ICEV} & 3.68 & 0.57 & 8.19 & 2.97 & 3.23 
\\ \bottomrule
\end{tabular}
\end{table}

We then computed the \ac{TCO} of the fleet  using existing retail prices presented in \Cref{mth_veh_features} (excluding value added tax, green incentives, etc.) and we plotted it in red in the lower part of \Cref{mth_group}. We note that without any additional taxes and incentives, it is not preferred to include any \acp{EV} in the fleet (comparing against the black line that indicates the \ac{TCO} of the fleet containing no \acp{EV}). 

We vary the price difference between the two competing classes of vehicles, namely  small \acp{ICEV} and \acp{EV}, by fixing the retail price of the former at 21k USD (where 1k = 1,000) and varying the price of the latter above 21k USD in \Cref{mth_group}. We note that \acp{EV} begin to become favorable only when the price difference drops below 5k USD. When the price difference falls below 3k USD, it is profitable to make all the small vehicles \acp{EV}. We note that significantly high savings are made by the firm by owning 44 \acp{EV} in the long run to the tune of 50k USD over owning 44 small \acp{ICEV}, when the price difference drops to 1000 USD.

\begin{figure}[ht!]
        \centering
        \begin{subfigure}[b]{0.5\textwidth}
                \centering                \includegraphics[scale = 0.5]{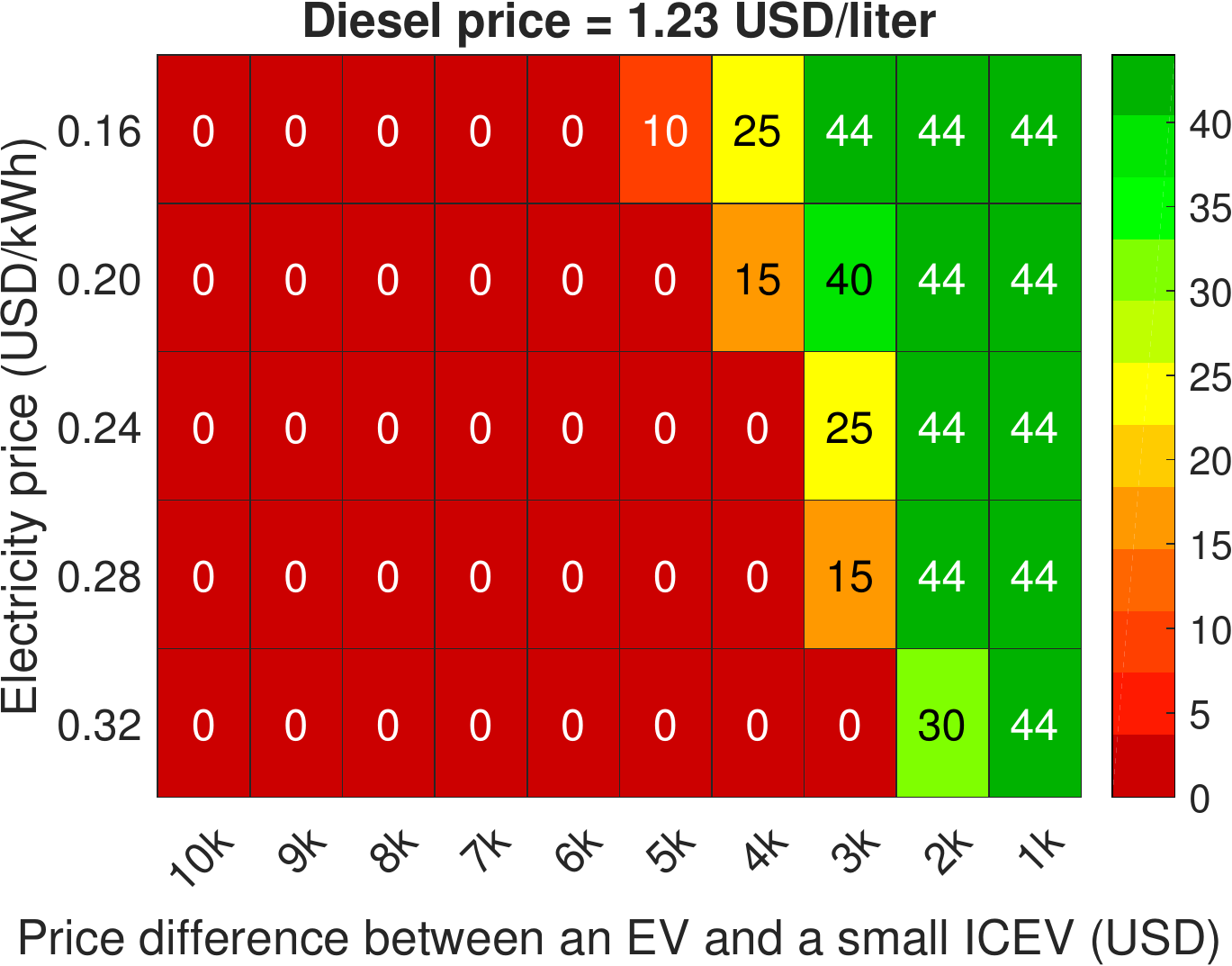}
                \caption{\small{1.23 USD}}
                \label{dsl_1.23}
        \end{subfigure}%
        \begin{subfigure}[b]{0.5\textwidth}
                \centering                \includegraphics[scale = 0.5]{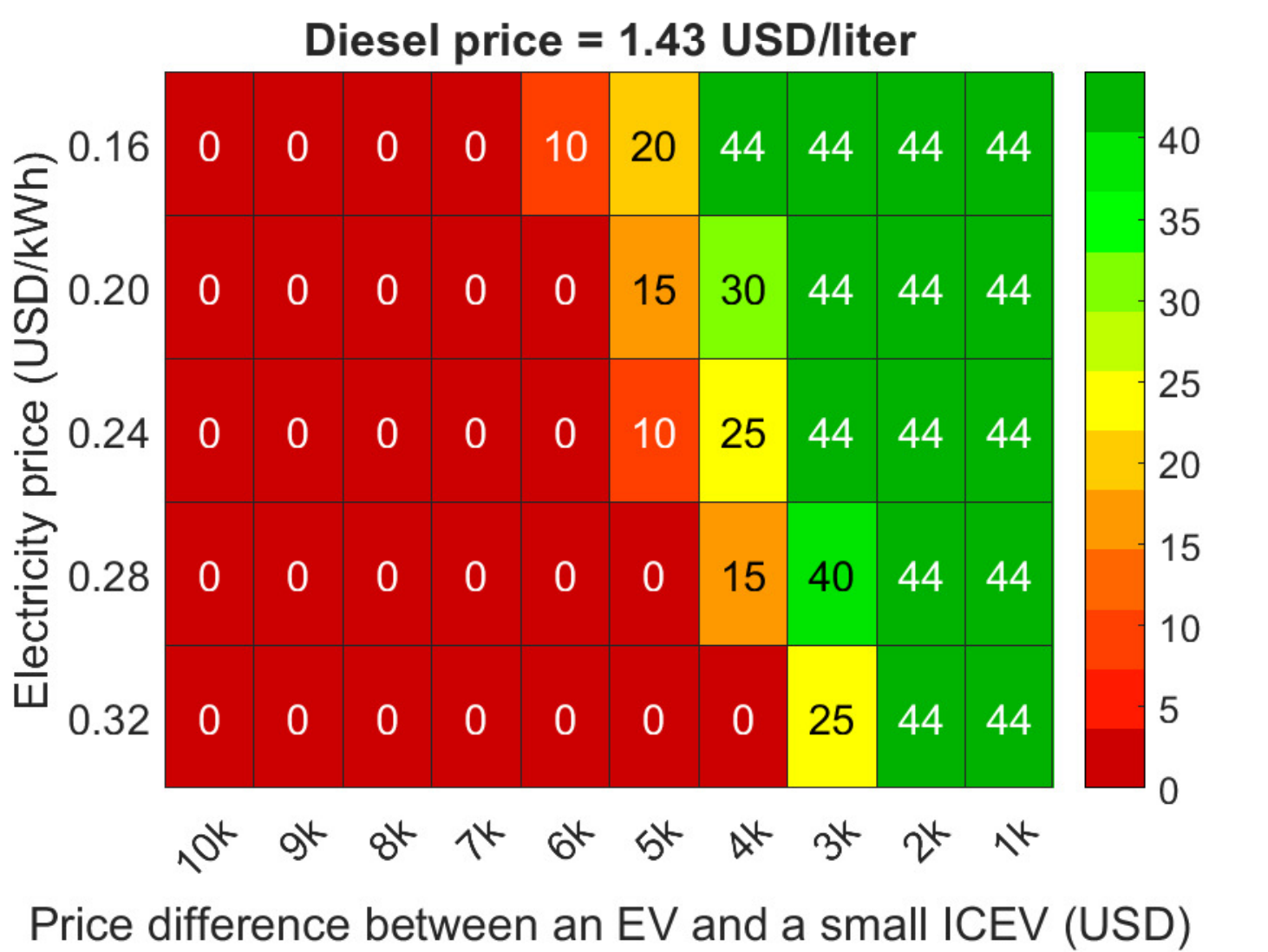}
                \caption{\small{1.43 USD}}
                \label{dsl_1.43}
        \end{subfigure}
        \begin{subfigure}[b]{0.5\textwidth}
                \centering                \includegraphics[scale = 0.5]{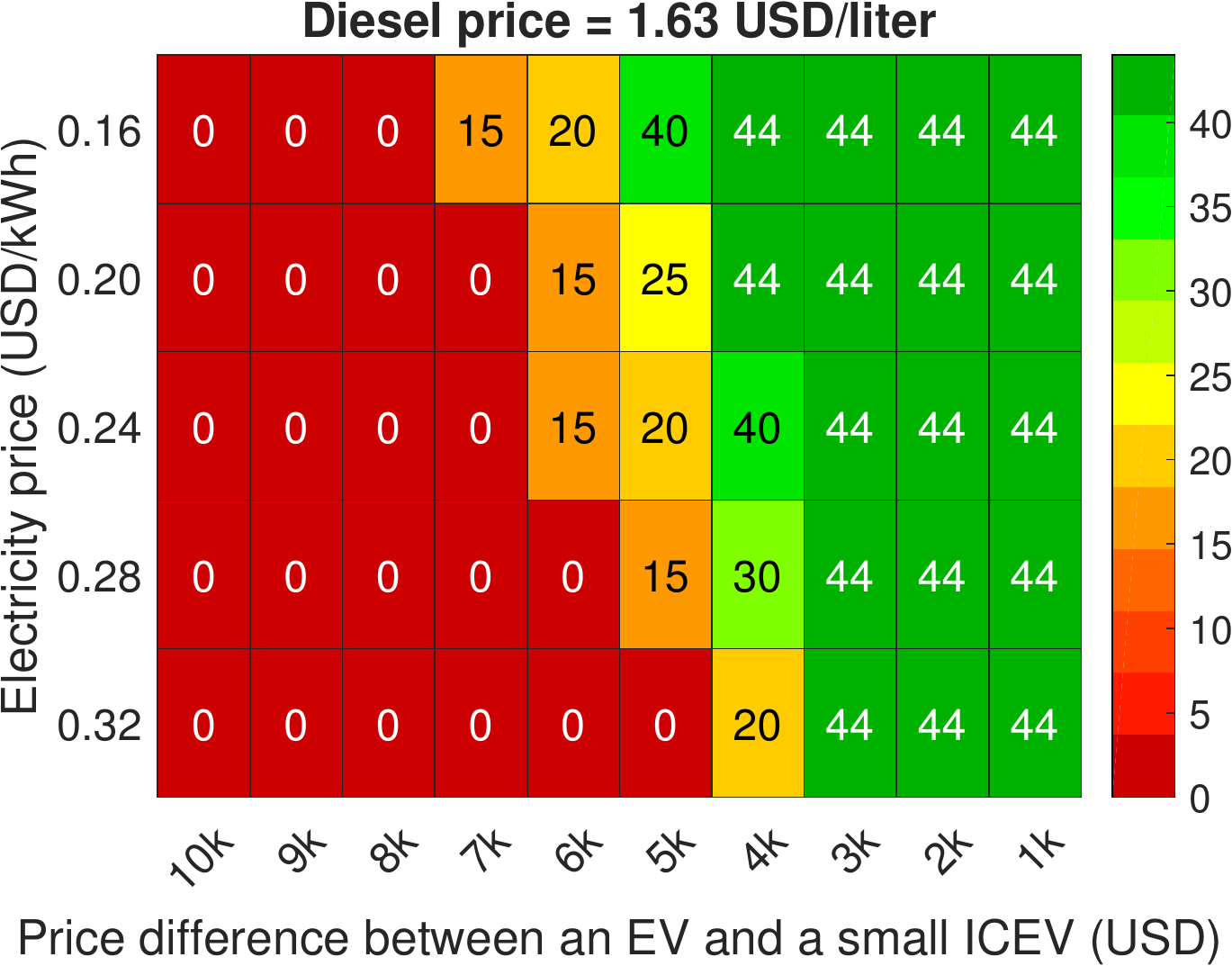}
                \caption{\small{1.63 USD}}
                \label{dsl_1.63}
        \end{subfigure}%
        \begin{subfigure}[b]{0.5\textwidth}
                \centering                \includegraphics[scale = 0.5]{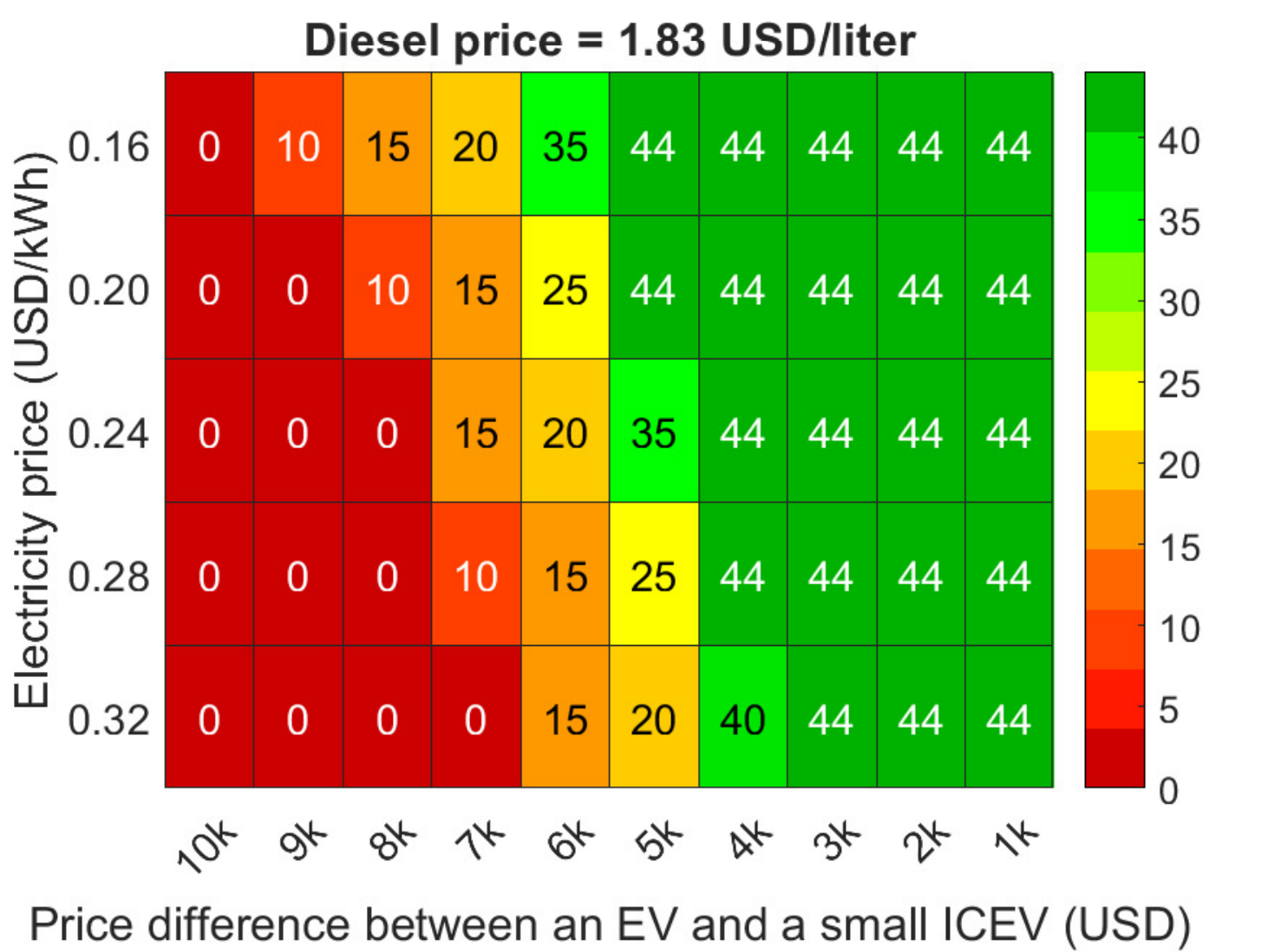}
                \caption{\small{1.83 USD}}
                \label{dsl_1.83}
        \end{subfigure}
        \caption{MTH case: Sensitivity of the tendency of \ac{EV} inclusion in the fleet to variation in electricity and diesel prices at various levels of price difference between \ac{EV} and small \ac{ICEV}: The number of \acp{EV} out of 44 small vehicles in the fleet up to which \ac{TCO} is lower than that of a no-\ac{EV} fleet mix is indicated. }
        \label{mth_sensitivity_analysis1}
\end{figure}

\begin{figure}[htbp!]
\begin{center}
\includegraphics[scale = 0.5]{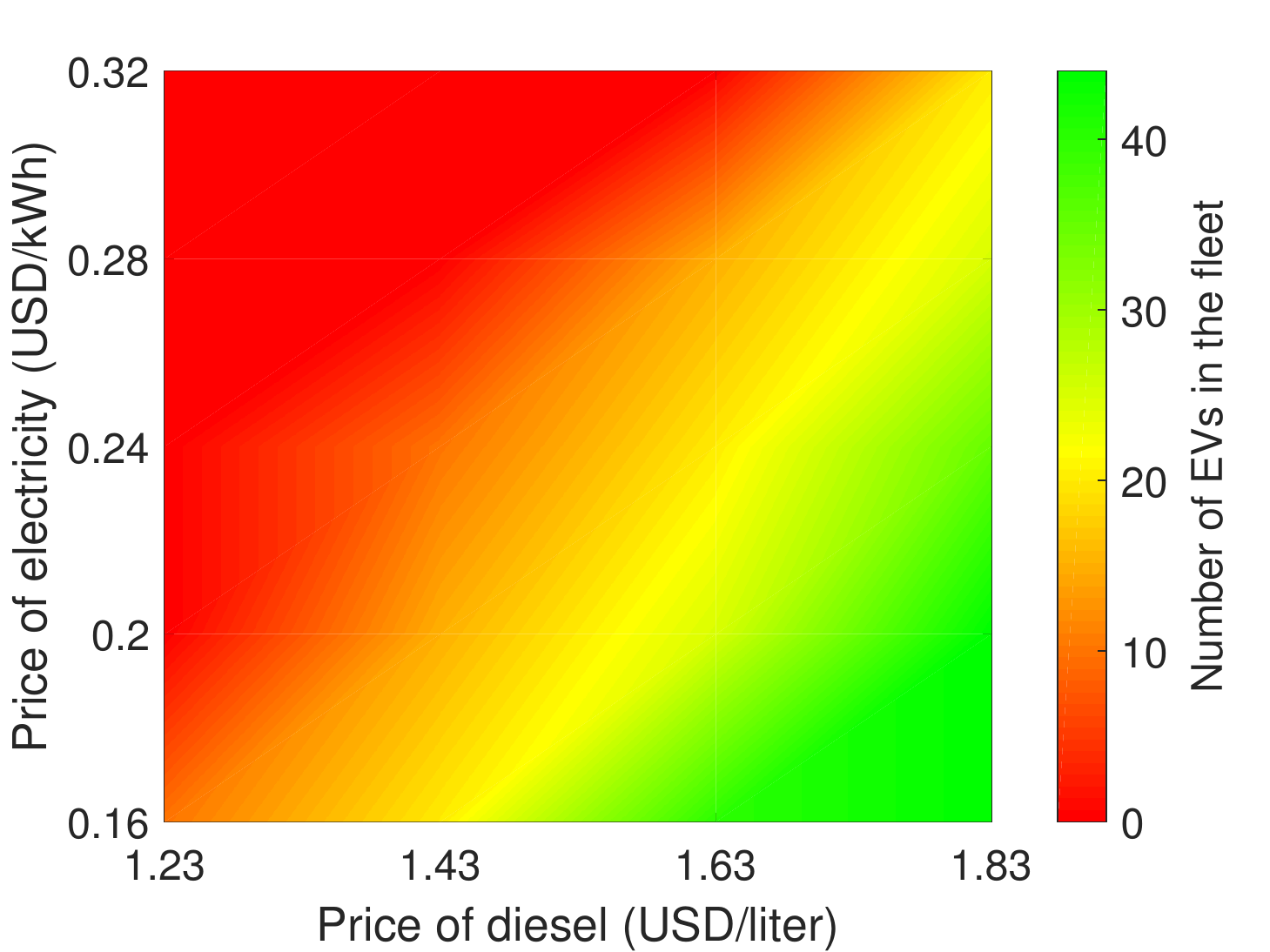}
\caption{MTH case: Variation of the number of \acp{EV} in the fleet in response to variation in diesel and electricity prices when \acp{EV} are 5k USD more expensive than \acp{ICEV}}
\label{mth_sensitivity_analysis2}
\end{center}
\end{figure}

Since the operational costs are heavily dependent on the prices of diesel and electricity, based on corresponding \acp{TCO}, we performed an analysis of the variation of the \acp{EV} composition in the fleet with variation in diesel prices between 1.23-1.83 USD/liter and electricity prices between 0.16-0.32 USD/kWh. \Cref{mth_sensitivity_analysis1} presents the number of \acp{EV} to be included in the fleet profitably at various levels of price difference between the competing vehicle classes. Confirming our intuition, we find that at lower electricity prices and higher diesel prices, there is an increased preference for the inclusion of \acp{EV} in the fleet. However, at the current retail price difference of 11k USD, it is not profitable to include \acp{EV} in the fleet for any reasonable combination of energy prices. Thus, additional incentives are necessary to encourage commercial \acp{EV} ownership. When the effective price difference is between 7k-3k USD, various partial \ac{EV} fleet mixes are optimal. Under a price difference of 2k USD,  all the small vehicles will be made \acp{EV} for almost all the energy price combinations. For a price difference of 5k USD, we present the variation of the tendency of electrification with the variation of energy source prices in \Cref{mth_sensitivity_analysis2}. As expected, lower electricity prices and higher diesel prices favor \ac{EV} strength in the fleet. 

\begin{figure}[ht!]
\begin{center}
\includegraphics[scale = 0.9]{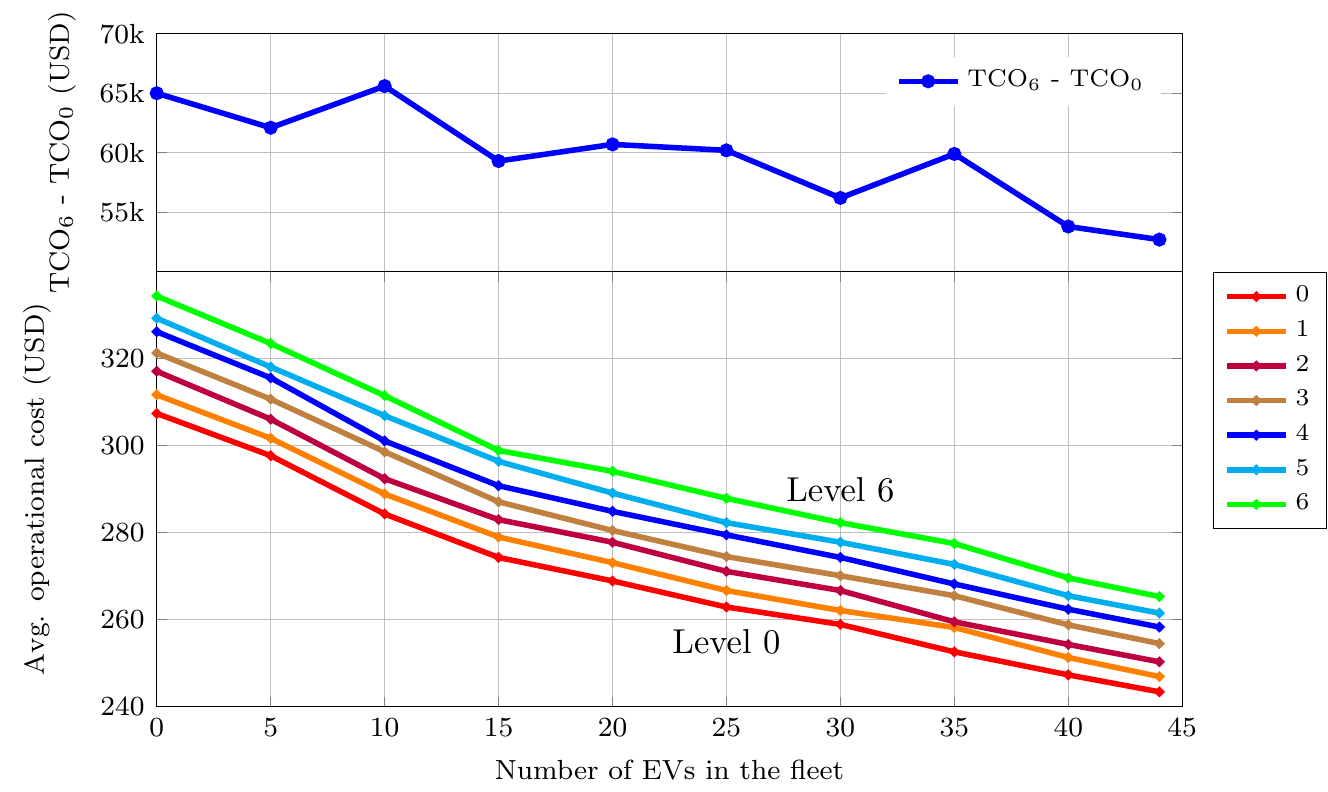}
\caption{MTH case: Variation of the average operational cost the fleet in response to variation of additional mass levels when temperature is fixed at 0 degree Celsius. In the upper plot, the maximum difference in \ac{TCO} for various fleet sizes is presented.  \ac{TCO}$_i$ is the \ac{TCO} computed at additional mass level $i$. }
\label{mth_exMass}
\end{center}
\end{figure}

Next, we vary the additional load carried in all the vehicles through seven levels $1, \dots, 6$, where the load mass is the level multiplied by 100 kg for \acp{EV} and small \acp{ICEV} while it is the level multiplied by 200 kg for the larger \acp{ICEV} (since the maximum allowed load for small and large vehicles is 600kg and  1200kg respectively). For these experiments, the ambient temperature is fixed at 0 degree Celsius.  From \Cref{mth_exMass}, we note that the average operational cost is affected by the additional mass level inasmuch as  leading to  9\% difference between the highest and lowest additional mass levels. The corresponding difference in \ac{TCO} is to the tune of 50-65k USD. Thus, we note that it is necessary to account for the additional mass carried by commercial vehicles while planning their operations.   

The solution method presented in this paper will serve as a useful tool to policy-makers in determining effective incentive schemes as well as to firms to determine the profitable extent of \ac{EV} ownership.

\section{Conclusion\label{conclusion}}
So far, fleet size and mix vehicle routing problems have not considered uncertainty of operational scenarios at the strategic decision stage. We fill this gap by proposing a novel two-stage stochastic program problem framework \ac{SFSMP} and solving it by \ac{SAA} through the use of the metaheuristic method \ac{ALNS} for operational cost determination and decision-making. From a modeling perspective, other recourse decisions, such as renting out acquired vehicles when idle or receiving chartered services when needed, may be investigated in future research. Also, the timing for selling acquired vehicles and the corresponding salvage values may be included in future work. Additionally, other solution methods may be explored to tackle the strategic as well as operational decision stages. 

In addition to presenting a novel problem, we evaluate the total cost of ownership of a mixed fleet with electric vehicles through two Danish case studies. The cases confirm the importance of considering the uncertainty in requests at the strategic planning stage and including cabin climate control power and auxiliary power in the computation of vehicle power for planning purposes. Additionally, the versatility of the proposed approach to not only determine the best fleet mix but also assess the robustness of the prescribed solution is demonstrated. We note that although we consider electric vehicles in this paper, other alternate sources of renewable energy (hydrogen-based, etc.) may be considered instead to understand their potential.
\section*{\normalsize{ACKNOWLEDGMENTS}}
The acknowledgements will be presented after peer review. 

\bibliography{references}

\begin{thebibliography}{73}
\expandafter\ifx\csname natexlab\endcsname\relax\def\natexlab#1{#1}\fi
\providecommand{\url}[1]{\texttt{#1}}
\providecommand{\href}[2]{#2}
\providecommand{\path}[1]{#1}
\providecommand{\DOIprefix}{doi:}
\providecommand{\ArXivprefix}{arXiv:}
\providecommand{\URLprefix}{URL: }
\providecommand{\Pubmedprefix}{pmid:}
\providecommand{\doi}[1]{\href{http://dx.doi.org/#1}{\path{#1}}}
\providecommand{\Pubmed}[1]{\href{pmid:#1}{\path{#1}}}
\providecommand{\bibinfo}[2]{#2}
\ifx\xfnm\relax \def\xfnm[#1]{\unskip,\space#1}\fi
\bibitem[{Ahmed and Shapiro(2002)}]{Ahmed2002TheRecourse}
\bibinfo{author}{Ahmed, S.}, \bibinfo{author}{Shapiro, A.},
  \bibinfo{year}{2002}.
\newblock \bibinfo{title}{{The sample average approximation method for
  stochastic programs with integer recourse}}.
\newblock \bibinfo{journal}{Technical Report, School of Industrial {\&} Systems
  Engineering, Georgia Institute of Technology.} , \bibinfo{pages}{1--24}.
\bibitem[{Anaya-Arenas et~al.(2014)Anaya-Arenas, Chabot, Renaud and
  Ruiz}]{Anaya-Arenas2014}
\bibinfo{author}{Anaya-Arenas, A.M.}, \bibinfo{author}{Chabot, T.},
  \bibinfo{author}{Renaud, J.}, \bibinfo{author}{Ruiz, A.},
  \bibinfo{year}{2014}.
\newblock \bibinfo{title}{{Biomedical sample transportation : A case study
  based on Qu{\'{e}}bec ’ s healthcare supply chain}}.
\bibitem[{Bakkehaug et~al.(2014)Bakkehaug, Eidem, Fagerholt and
  Hvattum}]{Bakkehaug2014AShipping}
\bibinfo{author}{Bakkehaug, R.}, \bibinfo{author}{Eidem, E.S.},
  \bibinfo{author}{Fagerholt, K.}, \bibinfo{author}{Hvattum, L.M.},
  \bibinfo{year}{2014}.
\newblock \bibinfo{title}{{A stochastic programming formulation for strategic
  fleet renewal in shipping}}.
\newblock \bibinfo{journal}{Transportation Research Part E: Logistics and
  Transportation Review} \bibinfo{volume}{72}, \bibinfo{pages}{60--76}.
\newblock \URLprefix \url{http://dx.doi.org/10.1016/j.tre.2014.09.010},
  \DOIprefix\doi{10.1016/j.tre.2014.09.010}.
\bibitem[{Baldacci et~al.(2009)Baldacci, Battarra and
  Vigo}]{Baldacci2009ValidCosts}
\bibinfo{author}{Baldacci, R.}, \bibinfo{author}{Battarra, M.},
  \bibinfo{author}{Vigo, D.}, \bibinfo{year}{2009}.
\newblock \bibinfo{title}{{Valid inequalities for the fleet size and mix
  vehicle routing problem with fixed costs}}.
\newblock \bibinfo{journal}{Networks} \bibinfo{volume}{54},
  \bibinfo{pages}{178--189}.
\bibitem[{Bekta{\c{s}} et~al.(2016)Bekta{\c{s}}, Demir and
  Laporte}]{Bektas2016GreenRouting}
\bibinfo{author}{Bekta{\c{s}}, T.}, \bibinfo{author}{Demir, E.},
  \bibinfo{author}{Laporte, G.}, \bibinfo{year}{2016}.
\newblock \bibinfo{title}{{Green vehicle routing}}, in:
  \bibinfo{booktitle}{Green Transportation Logistics: The Quest for Win-Win
  Solutions}. \bibinfo{publisher}{Springer}, pp. \bibinfo{pages}{243--265}.
\newblock \URLprefix
  \url{http://link.springer.com/10.1007/978-3-319-17175-3_7},
  \DOIprefix\doi{10.1007/978-3-319-17175-3{\_}7}.
\bibitem[{Bekta{\c{s}} et~al.(2018)Bekta{\c{s}}, Fabian~Ehmke, Psaraftis and
  Puchinger}]{Bektas2018TheTransportation}
\bibinfo{author}{Bekta{\c{s}}, T.}, \bibinfo{author}{Fabian~Ehmke, J.},
  \bibinfo{author}{Psaraftis, H.N.}, \bibinfo{author}{Puchinger, J.},
  \bibinfo{year}{2018}.
\newblock \bibinfo{title}{{The role of operational research in green freight
  transportation}}.
\newblock \bibinfo{journal}{European Journal of Operational Research}
  \bibinfo{volume}{16}, \bibinfo{pages}{1--17}.
\newblock \URLprefix \url{https://doi.org/10.1016/j.ejor.2018.06.001},
  \DOIprefix\doi{10.1016/j.ejor.2018.06.001}.
\bibitem[{Bekta{\c{s}} and Laporte(2011)}]{Bektas2011TheProblem}
\bibinfo{author}{Bekta{\c{s}}, T.}, \bibinfo{author}{Laporte, G.},
  \bibinfo{year}{2011}.
\newblock \bibinfo{title}{{The pollution-routing problem}}.
\newblock \bibinfo{journal}{Transportation Research Part B}
  \bibinfo{volume}{45}, \bibinfo{pages}{1232--1250}.
\newblock \URLprefix
  \url{https://ac.els-cdn.com/S019126151100018X/1-s2.0-S019126151100018X-main.pdf?_tid=4c1e55ee-1b2c-41cd-b056-79741a692377&acdnat=1538234693_52ac41d6ff6304b537ac7d15f4b955f9},
  \DOIprefix\doi{10.1016/j.trb.2011.02.004}.
\bibitem[{Christensen(2014)}]{Christensen2014AdaptiveThesis}
\bibinfo{author}{Christensen, J.M.}, \bibinfo{year}{2014}.
\newblock \bibinfo{title}{{Adaptive large neighbourhood search with exact
  methods for VRPTW (MSc Thesis)}}.
\newblock \URLprefix
  \url{http://production.datastore.cvt.dk/oafilestore?oid=575e9f9ee6f951534a00973b&targetid=56d754c0bf19455102000d9f}.
\bibitem[{Christiaens and Vanden~Berghe(2018)}]{ChristiaensJanSIbS}
\bibinfo{author}{Christiaens, J.}, \bibinfo{author}{Vanden~Berghe, G.},
  \bibinfo{year}{2018}.
\newblock \bibinfo{title}{{Slack induction by string removals for vehicle
  routing problems}}.
\newblock \bibinfo{type}{Technical Report}.
\newblock \URLprefix \url{[freely available]}.
\bibitem[{{CIVITAS}(2015)}]{CIVITAS2015MakingSustainable}
\bibinfo{author}{{CIVITAS}}, \bibinfo{year}{2015}.
\newblock \bibinfo{title}{{Making urban freight logistics more sustainable}}.
\newblock \bibinfo{journal}{CIVITAS Policy Note} ,
  \bibinfo{pages}{1--63}\URLprefix
  \url{http://www.eltis.org/resources/tools/civitas-policy-note-making-urban-freight-logistics-more-sustainable}.
\bibitem[{Conrad and Figliozzi(2011)}]{Conrad2011TheProblem}
\bibinfo{author}{Conrad, R.G.}, \bibinfo{author}{Figliozzi, M.A.},
  \bibinfo{year}{2011}.
\newblock \bibinfo{title}{{The recharging vehicle routing problem}}, in:
  \bibinfo{booktitle}{IIE Annual Conference. Proceedings; Norcross}, pp.
  \bibinfo{pages}{1--8}.
\newblock \URLprefix
  \url{http://search.proquest.com/docview/1190410233/abstract/F1BF6E5728C047DCPQ/1%0Ahttp://media.proquest.com/media/pq/classic/doc/2823141851/fmt/pi/rep/NONE?cit%3Aauth=Conrad%2C+Ryan+G%3BFigliozzi%2C+Miguel+Andres&cit%3Atitle=The+Recharging+Vehicle+Routing+Pro}.
\bibitem[{Cort{\'{e}}s et~al.(2011)Cort{\'{e}}s, Gendreau, Leng and
  Weintraub}]{Cortes2011ADemand}
\bibinfo{author}{Cort{\'{e}}s, C.E.}, \bibinfo{author}{Gendreau, M.},
  \bibinfo{author}{Leng, D.}, \bibinfo{author}{Weintraub, A.},
  \bibinfo{year}{2011}.
\newblock \bibinfo{title}{{A simulation-based approach for fleet design in a
  technician dispatch problem with stochastic demand}}.
\newblock \bibinfo{journal}{Journal of the Operational Research Society}
  \bibinfo{volume}{62}, \bibinfo{pages}{1510--1523}.
\newblock \DOIprefix\doi{10.1057/jors.2010.98}.
\bibitem[{Cort{\'{e}}s-Murcia et~al.(2019)Cort{\'{e}}s-Murcia, Prodhon and
  Murat~Afsar}]{Cortes-Murcia2019TheCustomers}
\bibinfo{author}{Cort{\'{e}}s-Murcia, D.L.}, \bibinfo{author}{Prodhon, C.},
  \bibinfo{author}{Murat~Afsar, H.}, \bibinfo{year}{2019}.
\newblock \bibinfo{title}{{The electric vehicle routing problem with time
  windows, partial recharges and satellite customers}}.
\newblock \bibinfo{journal}{Transportation Research Part E: Logistics and
  Transportation Review} \bibinfo{volume}{130}, \bibinfo{pages}{184--206}.
\newblock \URLprefix \url{https://doi.org/10.1016/j.tre.2019.08.015},
  \DOIprefix\doi{10.1016/j.tre.2019.08.015}.
\bibitem[{Crainic et~al.(2015)Crainic, Errico, Rei and
  Ricciardi}]{Crainic2015ModelingPlanning}
\bibinfo{author}{Crainic, T.G.}, \bibinfo{author}{Errico, F.},
  \bibinfo{author}{Rei, W.}, \bibinfo{author}{Ricciardi, N.},
  \bibinfo{year}{2015}.
\newblock \bibinfo{title}{{Modeling demand uncertainty in two-tier city
  logistics tactical planning}}.
\newblock \bibinfo{journal}{Transportation Science} \bibinfo{volume}{50},
  \bibinfo{pages}{559--578}.
\newblock \DOIprefix\doi{10.1287/trsc.2015.0606}.
\bibitem[{Crainic et~al.(2009)Crainic, Ricciardi and
  Storchi}]{Crainic2009ModelsSystems}
\bibinfo{author}{Crainic, T.G.}, \bibinfo{author}{Ricciardi, N.},
  \bibinfo{author}{Storchi, G.}, \bibinfo{year}{2009}.
\newblock \bibinfo{title}{{Models for evaluating and planning city logistics
  systems}}.
\newblock \bibinfo{journal}{Transportation Science} \bibinfo{volume}{43},
  \bibinfo{pages}{432--454}.
\newblock \DOIprefix\doi{10.1287/trsc.1090.0279}.
\bibitem[{Dell'Amico et~al.(2007)Dell'Amico, Monaci, Pagani and
  Vigo}]{DellAmico2007HeuristicWindows}
\bibinfo{author}{Dell'Amico, M.}, \bibinfo{author}{Monaci, M.},
  \bibinfo{author}{Pagani, C.}, \bibinfo{author}{Vigo, D.},
  \bibinfo{year}{2007}.
\newblock \bibinfo{title}{{Heuristic approaches for the fleet size and mix
  vehicle routing problem with time windows}}.
\newblock \bibinfo{journal}{Transportation Science} \bibinfo{volume}{41},
  \bibinfo{pages}{516--526}.
\newblock \URLprefix
  \url{http://pubsonline.informs.org/doi/abs/10.1287/trsc.1070.0190},
  \DOIprefix\doi{10.1287/trsc.1070.0190}.
\bibitem[{Demir et~al.(2011)Demir, Bektas and
  Laporte}]{Demir2011ATransportation}
\bibinfo{author}{Demir, E.}, \bibinfo{author}{Bektas, T.},
  \bibinfo{author}{Laporte, G.}, \bibinfo{year}{2011}.
\newblock \bibinfo{title}{{A comparative analysis of several vehicle emission
  models for road freight transportation}}.
\newblock \bibinfo{journal}{Transportation Research Part D}
  \bibinfo{volume}{16}, \bibinfo{pages}{347--357}.
\newblock \URLprefix
  \url{https://ac.els-cdn.com/S136192091100023X/1-s2.0-S136192091100023X-main.pdf?_tid=9a73b63d-4404-453c-b8d1-8740f1b072a7&acdnat=1544964742_cb13c82b7a30e2d53f6d3c4c6671b6ac},
  \DOIprefix\doi{10.1016/j.trd.2011.01.011}.
\bibitem[{Demir et~al.(2012)Demir, Bekta{\c{s}} and
  Laporte}]{Demir2012AnProblem}
\bibinfo{author}{Demir, E.}, \bibinfo{author}{Bekta{\c{s}}, T.},
  \bibinfo{author}{Laporte, G.}, \bibinfo{year}{2012}.
\newblock \bibinfo{title}{{An adaptive large neighborhood search heuristic for
  the pollution-routing problem}}.
\newblock \bibinfo{journal}{European Journal of Operational Research}
  \bibinfo{volume}{223}, \bibinfo{pages}{346--359}.
\newblock \DOIprefix\doi{10.1016/j.ejor.2012.06.044}.
\bibitem[{Dempster et~al.(1981)Dempster, Fisher, Jansen, Lageweg, Lenstra and
  Rinnooy~Kan}]{Dempster1981AnalyticalSystems}
\bibinfo{author}{Dempster, M.A.H.}, \bibinfo{author}{Fisher, M.L.},
  \bibinfo{author}{Jansen, L.}, \bibinfo{author}{Lageweg, B.J.},
  \bibinfo{author}{Lenstra, J.K.}, \bibinfo{author}{Rinnooy~Kan, A.H.G.},
  \bibinfo{year}{1981}.
\newblock \bibinfo{title}{{Analytical evaluation of hierarchical planning
  systems}}.
\newblock \bibinfo{journal}{Operations Research} \bibinfo{volume}{29},
  \bibinfo{pages}{707--716}.
\newblock \DOIprefix\doi{10.1287/opre.29.4.707}.
\bibitem[{Desrochers and Verhoog(1991)}]{Desrochers1991AProblem}
\bibinfo{author}{Desrochers, M.}, \bibinfo{author}{Verhoog, T.W.},
  \bibinfo{year}{1991}.
\newblock \bibinfo{title}{{A new heuristic for the fleet size and mix vehicle
  routing problem}}.
\newblock \bibinfo{journal}{Computers and Operations Research}
  \bibinfo{volume}{18}, \bibinfo{pages}{263--274}.
\newblock \DOIprefix\doi{10.1016/0305-0548(91)90028-P}.
\bibitem[{Dullaert et~al.(2002)Dullaert, Janssens, S{\"{o}}rensen and
  Vernimmen}]{Dullaert2002NewWindows}
\bibinfo{author}{Dullaert, W.}, \bibinfo{author}{Janssens, G.K.},
  \bibinfo{author}{S{\"{o}}rensen, K.}, \bibinfo{author}{Vernimmen, B.},
  \bibinfo{year}{2002}.
\newblock \bibinfo{title}{{New heuristics for the fleet size and mix vehicle
  routing problem with time windows}}.
\newblock \bibinfo{journal}{Journal of the Operational Research Society}
  \bibinfo{volume}{53}, \bibinfo{pages}{1232--1238}.
\bibitem[{Erdo{\u{g}}an and Miller-Hooks(2012)}]{Erdogan2012AProblem}
\bibinfo{author}{Erdo{\u{g}}an, S.}, \bibinfo{author}{Miller-Hooks, E.},
  \bibinfo{year}{2012}.
\newblock \bibinfo{title}{{A green vehicle routing problem}}.
\newblock \bibinfo{journal}{Transportation Research Part E}
  \bibinfo{volume}{48}, \bibinfo{pages}{100--114}.
\newblock \URLprefix
  \url{https://ac.els-cdn.com/S1366554511001062/1-s2.0-S1366554511001062-main.pdf?_tid=73cb4a45-2d42-42e9-9182-23626618bcc3&acdnat=1538061326_91d1be8958efc63f493cdffdfddabfdb},
  \DOIprefix\doi{10.1016/j.tre.2011.08.001}.
\bibitem[{{EUNADICS-AV}(2018)}]{EUNADICS-AV2018SolarSZA}
\bibinfo{author}{{EUNADICS-AV}}, \bibinfo{year}{2018}.
\newblock \bibinfo{title}{{Solar zenith angle (SZA)}}.
\newblock \URLprefix \url{http://sacs.aeronomie.be/info/sza.php}.
\bibitem[{Fayazbakhsh and Bahrami(2013)}]{Fayazbakhsh2013ComprehensiveMethod}
\bibinfo{author}{Fayazbakhsh, M.A.}, \bibinfo{author}{Bahrami, M.},
  \bibinfo{year}{2013}.
\newblock \bibinfo{title}{{Comprehensive modeling of vehicle air conditioning
  loads using heat balance method}}.
\newblock \bibinfo{journal}{SAE International} \URLprefix
  \url{http://papers.sae.org/2013-01-1507/},
  \DOIprefix\doi{10.4271/2013-01-1507}.
\bibitem[{Gendreau et~al.(2008)Gendreau, Hertz and
  Laporte}]{Gendreau2008AProblem}
\bibinfo{author}{Gendreau, M.}, \bibinfo{author}{Hertz, A.},
  \bibinfo{author}{Laporte, G.}, \bibinfo{year}{2008}.
\newblock \bibinfo{title}{{A tabu search heuristic for the vehicle routing
  problem}}.
\newblock \bibinfo{journal}{Management Science} \bibinfo{volume}{40},
  \bibinfo{pages}{1276--1290}.
\newblock \DOIprefix\doi{10.1287/mnsc.40.10.1276}.
\bibitem[{Gendreau et~al.(1999)Gendreau, Laporte, Musaraganyi and
  Taillard}]{Gendreau1999AProblem}
\bibinfo{author}{Gendreau, M.}, \bibinfo{author}{Laporte, G.},
  \bibinfo{author}{Musaraganyi, C.}, \bibinfo{author}{Taillard, E.D.},
  \bibinfo{year}{1999}.
\newblock \bibinfo{title}{{A tabu search heuristic for the heterogeneous fleet
  vehicle routing problem}}.
\newblock \bibinfo{journal}{Computers and Operations Research}
  \bibinfo{volume}{26}, \bibinfo{pages}{1153--1173}.
\newblock \DOIprefix\doi{10.1016/S0305-0548(98)00100-2}.
\bibitem[{Gheysens et~al.(1984)Gheysens, Golden and
  Assad}]{Gheysens1984AProblems}
\bibinfo{author}{Gheysens, F.}, \bibinfo{author}{Golden, B.},
  \bibinfo{author}{Assad, A.}, \bibinfo{year}{1984}.
\newblock \bibinfo{title}{{A comparison of techniques for solving the fleet
  size and mix vehicle routing problems}}.
\newblock \bibinfo{journal}{OR Spektrum} \bibinfo{volume}{6},
  \bibinfo{pages}{207--216}.
\newblock \DOIprefix\doi{10.1007/BF01720070}.
\bibitem[{Goeke and Schneider(2015)}]{Goeke2015RoutingVehicles}
\bibinfo{author}{Goeke, D.}, \bibinfo{author}{Schneider, M.},
  \bibinfo{year}{2015}.
\newblock \bibinfo{title}{{Routing a mixed fleet of electric and conventional
  vehicles}}.
\newblock \bibinfo{journal}{European Journal of Operational Research}
  \bibinfo{volume}{245}, \bibinfo{pages}{81--99}.
\newblock \URLprefix \url{http://dx.doi.org/10.1016/j.ejor.2015.01.049},
  \DOIprefix\doi{10.1016/j.ejor.2015.01.049}.
\bibitem[{Golden et~al.(1984)Golden, Assad, Levy and
  Gheysens}]{Golden1984TheProblem}
\bibinfo{author}{Golden, B.}, \bibinfo{author}{Assad, A.},
  \bibinfo{author}{Levy, L.}, \bibinfo{author}{Gheysens, F.},
  \bibinfo{year}{1984}.
\newblock \bibinfo{title}{{The fleet size and mix vehicle routing problem}}.
\newblock \bibinfo{journal}{Computers and Operations Research}
  \bibinfo{volume}{11}, \bibinfo{pages}{49--66}.
\newblock \DOIprefix\doi{10.1016/0305-0548(84)90007-8}.
\bibitem[{Grasas et~al.(2014)Grasas, Ramalhinho, Pessoa, Resende, Caball{\'{e}}
  and Barba}]{Grasas2014}
\bibinfo{author}{Grasas, A.}, \bibinfo{author}{Ramalhinho, H.},
  \bibinfo{author}{Pessoa, L.S.}, \bibinfo{author}{Resende, M.G.},
  \bibinfo{author}{Caball{\'{e}}, I.}, \bibinfo{author}{Barba, N.},
  \bibinfo{year}{2014}.
\newblock \bibinfo{title}{{On the improvement of blood sample collection at
  clinical laboratories}}.
\newblock \bibinfo{journal}{BMC Health Services Research} \bibinfo{volume}{14},
  \bibinfo{pages}{12}.
\newblock \URLprefix
  \url{http://bmchealthservres.biomedcentral.com/articles/10.1186/1472-6963-14-12},
  \DOIprefix\doi{10.1186/1472-6963-14-12}.
\bibitem[{Heinrich and Wichmann(2004)}]{Heinrich2004TrafficDisease}
\bibinfo{author}{Heinrich, J.}, \bibinfo{author}{Wichmann, H.E.},
  \bibinfo{year}{2004}.
\newblock \bibinfo{title}{{Traffic related pollutants in Europe and their
  effect on allergic disease}}.
\newblock \bibinfo{journal}{Current Opinion in Allergy and Clinical Immunology}
  \bibinfo{volume}{4}, \bibinfo{pages}{341--348}.
\newblock \URLprefix
  \url{http://www.embase.com/search/results?subaction=viewrecord&from=export&id=L39274073%5Cnhttp://dx.doi.org/10.1097/00130832-200410000-00003}.
\bibitem[{Hiermann et~al.(2019)Hiermann, Hartl, Puchinger and
  Vidal}]{Hiermann2019RoutingVehicles}
\bibinfo{author}{Hiermann, G.}, \bibinfo{author}{Hartl, R.F.},
  \bibinfo{author}{Puchinger, J.}, \bibinfo{author}{Vidal, T.},
  \bibinfo{year}{2019}.
\newblock \bibinfo{title}{{Routing a mix of conventional, plug-in hybrid, and
  electric vehicles}}.
\newblock \bibinfo{journal}{European Journal of Operational Research}
  \bibinfo{volume}{272}, \bibinfo{pages}{235--248}.
\newblock \URLprefix \url{https://doi.org/10.1016/j.ejor.2018.06.025},
  \DOIprefix\doi{10.1016/j.ejor.2018.06.025}.
\bibitem[{Hiermann et~al.(2016)Hiermann, Puchinger, Ropke and
  Hartl}]{Hiermann2016}
\bibinfo{author}{Hiermann, G.}, \bibinfo{author}{Puchinger, J.},
  \bibinfo{author}{Ropke, S.}, \bibinfo{author}{Hartl, R.F.},
  \bibinfo{year}{2016}.
\newblock \bibinfo{title}{{The electric fleet size and mix vehicle routing
  problem with time windows and recharging stations}}.
\newblock \bibinfo{journal}{European Journal of Operational Research}
  \bibinfo{volume}{252}, \bibinfo{pages}{995--1018}.
\newblock \URLprefix
  \url{https://www.sciencedirect.com/science/article/pii/S0377221716000837},
  \DOIprefix\doi{10.1016/J.EJOR.2016.01.038}.
\bibitem[{{International Renewable Energy
  Agency}(2017)}]{InternationalRenewableEnergyAgency2017ElectricBrief}
\bibinfo{author}{{International Renewable Energy Agency}},
  \bibinfo{year}{2017}.
\newblock \bibinfo{title}{{Electric vehicles: Technology brief}}.
\newblock \bibinfo{type}{Technical Report}.
\newblock \URLprefix \url{www.irena.org}.
\bibitem[{Kergosien et~al.(2014)Kergosien, Ruiz and
  Soriano}]{Kergosien2014AServices}
\bibinfo{author}{Kergosien, Y.}, \bibinfo{author}{Ruiz, A.},
  \bibinfo{author}{Soriano, P.}, \bibinfo{year}{2014}.
\newblock \bibinfo{title}{{A routing problem for medical test sample collection
  in home health care services}}, in: \bibinfo{booktitle}{Proceedings of the
  International Conference on Health Care Systems Engineering}.
\newblock \DOIprefix\doi{10.1007/978-3-319-01848-5}.
\bibitem[{Keskin and {\c{C}}atay(2016)}]{Keskin2016PartialWindows}
\bibinfo{author}{Keskin, M.}, \bibinfo{author}{{\c{C}}atay, B.},
  \bibinfo{year}{2016}.
\newblock \bibinfo{title}{{Partial recharge strategies for the electric vehicle
  routing problem with time windows}}.
\newblock \bibinfo{journal}{Transportation Research Part C}
  \bibinfo{volume}{65}, \bibinfo{pages}{111--127}.
\newblock \URLprefix \url{http://dx.doi.org/10.1016/j.trc.2016.01.013},
  \DOIprefix\doi{10.1016/j.trc.2016.01.013}.
\bibitem[{Kijewska et~al.(2016)Kijewska, Konicki and
  Iwan}]{Kijewska2016FreightExample}
\bibinfo{author}{Kijewska, K.}, \bibinfo{author}{Konicki, W.},
  \bibinfo{author}{Iwan, S.}, \bibinfo{year}{2016}.
\newblock \bibinfo{title}{{Freight transport pollution propagation at urban
  areas based on Szczecin example}}.
\newblock \bibinfo{journal}{Transportation Research Procedia}
  \bibinfo{volume}{14}, \bibinfo{pages}{1543--1552}.
\newblock \DOIprefix\doi{10.1016/j.trpro.2016.05.119}.
\bibitem[{K{\"{u}}hlwein(2016)}]{Kuhlwein2016DRIVING2025}
\bibinfo{author}{K{\"{u}}hlwein, J.}, \bibinfo{year}{2016}.
\newblock \bibinfo{title}{{Driving resistances of light-duty vehicles in
  Europe: present situation, trends, and scenarios for 2025}}.
\newblock \bibinfo{type}{Technical Report}.
\newblock \URLprefix \url{www.theicct.org}.
\bibitem[{Lau et~al.(2003)Lau, Sim and Teo}]{Lau2003VehicleVehicles}
\bibinfo{author}{Lau, H.C.}, \bibinfo{author}{Sim, M.}, \bibinfo{author}{Teo,
  K.M.}, \bibinfo{year}{2003}.
\newblock \bibinfo{title}{{Vehicle routing problem with time windows and a
  limited number of vehicles}}.
\newblock \bibinfo{journal}{European Journal of Operational Research}
  \bibinfo{volume}{148}, \bibinfo{pages}{559--569}.
\newblock \DOIprefix\doi{10.1016/S0377-2217(02)00363-6}.
\bibitem[{Lin et~al.(2016)Lin, Zhou and Wolfson}]{Lin2016}
\bibinfo{author}{Lin, J.}, \bibinfo{author}{Zhou, W.},
  \bibinfo{author}{Wolfson, O.}, \bibinfo{year}{2016}.
\newblock \bibinfo{title}{{Electric vehicle routing problem}}.
\newblock \bibinfo{journal}{Transportation Research Procedia}
  \bibinfo{volume}{12}, \bibinfo{pages}{508--521}.
\newblock \DOIprefix\doi{10.1016/j.trpro.2016.02.007}.
\bibitem[{Liu and Shen(1999)}]{Liu1999TheWindows}
\bibinfo{author}{Liu, F.H.}, \bibinfo{author}{Shen, S.Y.},
  \bibinfo{year}{1999}.
\newblock \bibinfo{title}{{The fleet size and mix vehicle routing problem with
  time windows}}.
\newblock \bibinfo{journal}{The Journal of the Operational Research Society}
  \bibinfo{volume}{50}, \bibinfo{pages}{721--732}.
\newblock \URLprefix
  \url{http://www.jstor.org/stable/3010326?seq=1#page_scan_tab_contents}.
\bibitem[{Liu et~al.(2018)Liu, Wang, Yamamoto and
  Morikawa}]{Liu2018ExploringConsumption}
\bibinfo{author}{Liu, K.}, \bibinfo{author}{Wang, J.},
  \bibinfo{author}{Yamamoto, T.}, \bibinfo{author}{Morikawa, T.},
  \bibinfo{year}{2018}.
\newblock \bibinfo{title}{{Exploring the interactive effects of ambient
  temperature and vehicle auxiliary loads on electric vehicle energy
  consumption}}.
\newblock \bibinfo{journal}{Applied Energy} \bibinfo{volume}{227},
  \bibinfo{pages}{324--331}.
\newblock \URLprefix \url{http://dx.doi.org/10.1016/j.apenergy.2017.08.074},
  \DOIprefix\doi{10.1016/j.apenergy.2017.08.074}.
\bibitem[{Luxen and Vetter(2011)}]{luxen-vetter-2011}
\bibinfo{author}{Luxen, D.}, \bibinfo{author}{Vetter, C.},
  \bibinfo{year}{2011}.
\newblock \bibinfo{title}{{Real-time routing with OpenStreetMap data}}, in:
  \bibinfo{booktitle}{Proceedings of the 19th ACM SIGSPATIAL International
  Conference on Advances in Geographic Information Systems},
  \bibinfo{publisher}{ACM}, \bibinfo{address}{New York, NY, USA}. pp.
  \bibinfo{pages}{513--516}.
\newblock \URLprefix \url{http://doi.acm.org/10.1145/2093973.2094062},
  \DOIprefix\doi{10.1145/2093973.2094062}.
\bibitem[{Mitropoulos et~al.(2017)Mitropoulos, Prevedouros and
  Kopelias}]{Mitropoulos2017TotalVehicle}
\bibinfo{author}{Mitropoulos, L.K.}, \bibinfo{author}{Prevedouros, P.D.},
  \bibinfo{author}{Kopelias, P.}, \bibinfo{year}{2017}.
\newblock \bibinfo{title}{{Total cost of ownership and externalities of
  conventional, hybrid and electric vehicle}}.
\newblock \bibinfo{journal}{Transportation Research Procedia}
  \bibinfo{volume}{24}, \bibinfo{pages}{267--274}.
\newblock \URLprefix \url{http://dx.doi.org/10.1016/j.trpro.2017.05.117},
  \DOIprefix\doi{10.1016/j.trpro.2017.05.117}.
\bibitem[{Montoya et~al.(2017)Montoya, Gu{\'{e}}ret, Mendoza and
  Villegas}]{Montoya2017}
\bibinfo{author}{Montoya, A.}, \bibinfo{author}{Gu{\'{e}}ret, C.},
  \bibinfo{author}{Mendoza, J.E.}, \bibinfo{author}{Villegas, J.G.},
  \bibinfo{year}{2017}.
\newblock \bibinfo{title}{{The electric vehicle routing problem with nonlinear
  charging function}}.
\newblock \bibinfo{journal}{Transportation Research Part B: Methodological}
  \DOIprefix\doi{10.1016/j.trb.2017.02.004}.
\bibitem[{Murakami(2017)}]{Murakami2017ARouting}
\bibinfo{author}{Murakami, K.}, \bibinfo{year}{2017}.
\newblock \bibinfo{title}{{A new model and approach to electric and
  diesel-powered vehicle routing}}.
\newblock \bibinfo{journal}{Transportation Research Part E}
  \bibinfo{volume}{107}, \bibinfo{pages}{23--37}.
\newblock \URLprefix
  \url{https://linkinghub.elsevier.com/retrieve/pii/S136655451730073X},
  \DOIprefix\doi{10.1016/j.tre.2017.09.004}.
\bibitem[{Osman and Salhi(1996)}]{Osman1996LocalProblem}
\bibinfo{author}{Osman, I.H.}, \bibinfo{author}{Salhi, S.},
  \bibinfo{year}{1996}.
\newblock \bibinfo{title}{{Local search strategies for the mix fleet routing
  problem}}.
\newblock \bibinfo{journal}{Modern Heuristics Search Methods} ,
  \bibinfo{pages}{131--153}.
\bibitem[{Pasha et~al.(2016)Pasha, Hoff and Hvattum}]{Pasha2016SimpleProblem}
\bibinfo{author}{Pasha, U.}, \bibinfo{author}{Hoff, A.},
  \bibinfo{author}{Hvattum, L.M.}, \bibinfo{year}{2016}.
\newblock \bibinfo{title}{{Simple heuristics for the multi-period fleet size
  and mix vehicle routing problem}}.
\newblock \bibinfo{journal}{INFOR} \bibinfo{volume}{54},
  \bibinfo{pages}{97--120}.
\newblock \URLprefix \url{http://dx.doi.org/10.1080/03155986.2016.1149314},
  \DOIprefix\doi{10.1080/03155986.2016.1149314}.
\bibitem[{Pelletier et~al.(2016)Pelletier, Jabali and Laporte}]{Pelletier2016}
\bibinfo{author}{Pelletier, S.}, \bibinfo{author}{Jabali, O.},
  \bibinfo{author}{Laporte, G.}, \bibinfo{year}{2016}.
\newblock \bibinfo{title}{{50th anniversary invited article - Goods
  distribution with electric vehicles: review and research perspectives}}.
\newblock \bibinfo{journal}{Transportation Science} \bibinfo{volume}{50},
  \bibinfo{pages}{3--22}.
\newblock \URLprefix
  \url{http://pubsonline.informs.org/doi/10.1287/trsc.2015.0646},
  \DOIprefix\doi{10.1287/trsc.2015.0646}.
\bibitem[{Pelletier et~al.(2017)Pelletier, Jabali, Laporte and
  Veneroni}]{Pelletier2017BatteryModels}
\bibinfo{author}{Pelletier, S.}, \bibinfo{author}{Jabali, O.},
  \bibinfo{author}{Laporte, G.}, \bibinfo{author}{Veneroni, M.},
  \bibinfo{year}{2017}.
\newblock \bibinfo{title}{{Battery degradation and behaviour for electric
  vehicles: Review and numerical analyses of several models}}.
\newblock \bibinfo{journal}{Transportation Research Part B: Methodological}
  \bibinfo{volume}{103}, \bibinfo{pages}{158--187}.
\newblock \DOIprefix\doi{10.1016/j.trb.2017.01.020}.
\bibitem[{Pillac et~al.(2013)Pillac, Gu{\'{e}}ret and
  Medaglia}]{Pillac2013AProblem}
\bibinfo{author}{Pillac, V.}, \bibinfo{author}{Gu{\'{e}}ret, C.},
  \bibinfo{author}{Medaglia, A.L.}, \bibinfo{year}{2013}.
\newblock \bibinfo{title}{{A parallel matheuristic for the technician routing
  and scheduling problem}}.
\newblock \bibinfo{journal}{Optimization Letters} \bibinfo{volume}{7},
  \bibinfo{pages}{1525--1535}.
\newblock \DOIprefix\doi{10.1007/s11590-012-0567-4}.
\bibitem[{Pisinger and Ropke(2007)}]{Pisinger2007AProblems}
\bibinfo{author}{Pisinger, D.}, \bibinfo{author}{Ropke, S.},
  \bibinfo{year}{2007}.
\newblock \bibinfo{title}{{A general heuristic for vehicle routing problems}}.
\newblock \bibinfo{journal}{Computers and Operations Research}
  \bibinfo{volume}{34}, \bibinfo{pages}{2403--2435}.
\newblock \DOIprefix\doi{10.1016/j.cor.2005.09.012}.
\bibitem[{Prins(2004)}]{Prins2004AProblem}
\bibinfo{author}{Prins, C.}, \bibinfo{year}{2004}.
\newblock \bibinfo{title}{{A simple and effective evolutionary algorithm for
  the vehicle routing problem}}.
\newblock \bibinfo{journal}{Computers and Operations Research}
  \bibinfo{volume}{31}, \bibinfo{pages}{1985--2002}.
\newblock \DOIprefix\doi{10.1016/S0305-0548(03)00158-8}.
\bibitem[{Ramsey and Kuehn(2018)}]{Ramsey2018SolarRadiation}
\bibinfo{author}{Ramsey, J.}, \bibinfo{author}{Kuehn, T.},
  \bibinfo{year}{2018}.
\newblock \bibinfo{title}{{Solar radiation}} , \bibinfo{pages}{1--14}\URLprefix
  \url{http://www.me.umn.edu/courses/me4131/LabManual/AppDSolarRadiation.pdf}.
\bibitem[{Renaud and Boctor(2002)}]{Renaud2002AProblem}
\bibinfo{author}{Renaud, J.}, \bibinfo{author}{Boctor, F.F.},
  \bibinfo{year}{2002}.
\newblock \bibinfo{title}{{A sweep-based algorithm for the fleet size and mix
  vehicle routing problem}}.
\newblock \bibinfo{journal}{European Journal of Operational Research}
  \bibinfo{volume}{140}, \bibinfo{pages}{618--628}.
\newblock \DOIprefix\doi{10.1016/S0377-2217(01)00237-5}.
\bibitem[{{Renault Danmark}(2018)}]{RenaultDanmark2018NyDanmark}
\bibinfo{author}{{Renault Danmark}}, \bibinfo{year}{2018}.
\newblock \bibinfo{title}{{Ny Kangoo Z.E. | Elbiler | Renault Danmark}}.
\newblock \URLprefix
  \url{https://www.renault.dk/biler/varebiler/ny-kangoo-ze.html}.
\bibitem[{Rich and Hansen(2016)}]{Rich2016TheResults}
\bibinfo{author}{Rich, J.}, \bibinfo{author}{Hansen, C.O.},
  \bibinfo{year}{2016}.
\newblock \bibinfo{title}{{The Danish national passenger model – Model
  specification and results}}.
\newblock \bibinfo{journal}{European Journal of Transport and Infrastructure
  Research} \bibinfo{volume}{16}, \bibinfo{pages}{573--599}.
\bibitem[{Rogge et~al.(2018)Rogge, van~der Hurk, Larsen and
  Sauer}]{Rogge2018ElectricInfrastructure}
\bibinfo{author}{Rogge, M.}, \bibinfo{author}{van~der Hurk, E.},
  \bibinfo{author}{Larsen, A.}, \bibinfo{author}{Sauer, D.U.},
  \bibinfo{year}{2018}.
\newblock \bibinfo{title}{{Electric bus fleet size and mix problem with
  optimization of charging infrastructure}}.
\newblock \bibinfo{journal}{Applied Energy} \bibinfo{volume}{211},
  \bibinfo{pages}{282--295}.
\newblock \URLprefix
  \url{https://linkinghub.elsevier.com/retrieve/pii/S0306261917316355},
  \DOIprefix\doi{10.1016/j.apenergy.2017.11.051}.
\bibitem[{R{\o}pke and Pisinger(2006)}]{Rpke2006AnWindows}
\bibinfo{author}{R{\o}pke, S.}, \bibinfo{author}{Pisinger, D.},
  \bibinfo{year}{2006}.
\newblock \bibinfo{title}{{An adaptive large neighborhood search heuristic for
  the pickup and delivery problem with time windows}}.
\newblock \bibinfo{journal}{Transportation Science} \URLprefix
  \url{http://orbit.dtu.dk/files/3154899/An adaptive large neighborhood search
  heuristic for the pickup and delivery problem with time
  windows_TechRep_ropke_pisinger.pdf}, \DOIprefix\doi{10.1287/trsc.1050.0135}.
\bibitem[{Savelsbergh and Van~Woensel(2016)}]{Savelsbergh201650thOpportunities}
\bibinfo{author}{Savelsbergh, M.}, \bibinfo{author}{Van~Woensel, T.},
  \bibinfo{year}{2016}.
\newblock \bibinfo{title}{{50th anniversary invited article - City logistics:
  Challenges and opportunities}}.
\newblock \bibinfo{journal}{Transportation Science} \bibinfo{volume}{50},
  \bibinfo{pages}{579--590}.
\newblock \URLprefix
  \url{http://pubsonline.informs.org/doi/10.1287/trsc.2016.0675},
  \DOIprefix\doi{10.1287/trsc.2016.0675}.
\bibitem[{Schiffer et~al.(2018)Schiffer, St{\"{u}}tz and
  Walther}]{Schiffer2018ElectricNetworks}
\bibinfo{author}{Schiffer, M.}, \bibinfo{author}{St{\"{u}}tz, S.},
  \bibinfo{author}{Walther, G.}, \bibinfo{year}{2018}.
\newblock \bibinfo{title}{{Electric commercial vehicles in mid-haul logistics
  networks}}.
\newblock \bibinfo{journal}{Green Energy and Technology} ,
  \bibinfo{pages}{153--173}\DOIprefix\doi{10.1007/978-3-319-69950-9{\_}7}.
\bibitem[{Schneider et~al.(2014)Schneider, Stenger and
  Goeke}]{Schneider2014TheStations}
\bibinfo{author}{Schneider, M.}, \bibinfo{author}{Stenger, A.},
  \bibinfo{author}{Goeke, D.}, \bibinfo{year}{2014}.
\newblock \bibinfo{title}{{The electric vehicle-routing problem with time
  windows and recharging stations}}.
\newblock \bibinfo{journal}{Transportation Science} \bibinfo{volume}{48},
  \bibinfo{pages}{500--520}.
\newblock \URLprefix \url{http://pubsonline.informs.org},
  \DOIprefix\doi{10.1287/trsc.2013.0490}.
\bibitem[{Schultz(1996)}]{Schultz1996RatesRecourse}
\bibinfo{author}{Schultz, R.}, \bibinfo{year}{1996}.
\newblock \bibinfo{title}{{Rates of convergence in stochastic programs with
  complete integer recourse}}.
\newblock \bibinfo{journal}{SIAM Journal on Optimization} \bibinfo{volume}{6},
  \bibinfo{pages}{1138--1152}.
\bibitem[{Shapiro et~al.(20009)Shapiro, Dentcheva and
  Ruszczy{\'{n}}ski}]{Shapiro20009LecturesProgramming}
\bibinfo{author}{Shapiro, A.}, \bibinfo{author}{Dentcheva, D.},
  \bibinfo{author}{Ruszczy{\'{n}}ski, A.}, \bibinfo{year}{20009}.
\newblock \bibinfo{title}{{Lectures on stochastic programming}}.
\newblock \bibinfo{number}{January 2009}.
\newblock \DOIprefix\doi{10.1137/1.9780898718751}.
\bibitem[{Song et~al.(2015)Song, Kwon and Kim}]{Song2015AirVehicle}
\bibinfo{author}{Song, B.}, \bibinfo{author}{Kwon, J.}, \bibinfo{author}{Kim,
  Y.}, \bibinfo{year}{2015}.
\newblock \bibinfo{title}{{Air conditioning system sizing for pure electric
  vehicle}}.
\newblock \bibinfo{journal}{World Electric Vehicle Journal}
  \bibinfo{volume}{7}, \bibinfo{pages}{407--413}.
\newblock \DOIprefix\doi{10.4271/960688}.
\bibitem[{St{\aa}lhane et~al.(2019)St{\aa}lhane, Halvorsen-Weare, Non{\aa}s and
  Pantuso}]{Stalhane2019OptimizingFarms}
\bibinfo{author}{St{\aa}lhane, M.}, \bibinfo{author}{Halvorsen-Weare, E.E.},
  \bibinfo{author}{Non{\aa}s, L.M.}, \bibinfo{author}{Pantuso, G.},
  \bibinfo{year}{2019}.
\newblock \bibinfo{title}{{Optimizing vessel fleet size and mix to support
  maintenance operations at offshore wind farms}}.
\newblock \bibinfo{journal}{European Journal of Operational Research}
  \bibinfo{volume}{276}, \bibinfo{pages}{495--509}.
\newblock \DOIprefix\doi{10.1016/j.ejor.2019.01.023}.
\bibitem[{Taillard(1999)}]{Taillard1999AVRP}
\bibinfo{author}{Taillard, E.}, \bibinfo{year}{1999}.
\newblock \bibinfo{title}{{A heuristic column generation method for the
  heterogeneous fleet VRP}}.
\newblock \bibinfo{journal}{VRP. RAIRO} \bibinfo{volume}{1},
  \bibinfo{pages}{1--14}.
\bibitem[{Tang and Veelenturf(2019)}]{Tang2019TheEra}
\bibinfo{author}{Tang, C.S.}, \bibinfo{author}{Veelenturf, L.P.},
  \bibinfo{year}{2019}.
\newblock \bibinfo{title}{{The strategic role of logistics in the industry 4.0
  era}}.
\newblock \bibinfo{journal}{Transportation Research Part E: Logistics and
  Transportation Review} \bibinfo{volume}{129}, \bibinfo{pages}{1--11}.
\newblock \URLprefix \url{https://doi.org/10.1016/j.tre.2019.06.004},
  \DOIprefix\doi{10.1016/j.tre.2019.06.004}.
\bibitem[{{UN}(2018)}]{UN201868Affairs}
\bibinfo{author}{{UN}}, \bibinfo{year}{2018}.
\newblock \bibinfo{title}{{68{\%} of the world population projected to live in
  urban areas by 2050, says UN | UN DESA | United Nations Department of
  Economic and Social Affairs}}.
\newblock \URLprefix
  \url{https://www.un.org/development/desa/en/news/population/2018-revision-of-world-urbanization-prospects.html}.
\bibitem[{Valentina et~al.(2014)Valentina, Viehl, Bringmann and
  Rosenstiel}]{Valentina2014HVACVehicles}
\bibinfo{author}{Valentina, R.}, \bibinfo{author}{Viehl, A.},
  \bibinfo{author}{Bringmann, O.}, \bibinfo{author}{Rosenstiel, W.},
  \bibinfo{year}{2014}.
\newblock \bibinfo{title}{{HVAC system modeling for range prediction of
  electric vehicles}}.
\newblock \bibinfo{journal}{IEEE Intelligent Vehicles Symposium, Proceedings} ,
  \bibinfo{pages}{1145--1150}\DOIprefix\doi{10.1109/IVS.2014.6856500}.
\bibitem[{Villegas et~al.(2018)Villegas, Gu{\'{e}}ret, Mendoza, Montoya and
  Villegas}]{Villegas2018TheVehicle}
\bibinfo{author}{Villegas, J.}, \bibinfo{author}{Gu{\'{e}}ret, C.},
  \bibinfo{author}{Mendoza, J.E.}, \bibinfo{author}{Montoya, A.},
  \bibinfo{author}{Villegas, J.G.}, \bibinfo{year}{2018}.
\newblock \bibinfo{title}{{The technician routing and scheduling problem with
  conventional and electric vehicle}}.
\newblock \bibinfo{type}{Technical Report}.
\newblock \URLprefix \url{https://hal.archives-ouvertes.fr/hal-01813887}.
\bibitem[{Waraich et~al.(2015)Waraich, Georges, Galus and
  Axhausen}]{Waraich2015}
\bibinfo{author}{Waraich, R.A.}, \bibinfo{author}{Georges, G.},
  \bibinfo{author}{Galus, M.D.}, \bibinfo{author}{Axhausen, K.W.},
  \bibinfo{year}{2015}.
\newblock \bibinfo{title}{{Adding electric vehicle modeling capability to an
  agent-based transport simulation}}, in: \bibinfo{booktitle}{Transportation
  Systems and Engineering: Concepts, Methodologies, Tools, and Applications},
  pp. \bibinfo{pages}{282--318}.
\newblock \URLprefix
  \url{https://www.scopus.com/inward/record.uri?eid=2-s2.0-84958914772&partnerID=40&md5=e9ba4c006fb925f5d3eec827b6c9cd11},
  \DOIprefix\doi{10.4018/978-1-4666-8473-7.ch075}.
\bibitem[{Zufferey et~al.(2016)Zufferey, Cho and Glardon}]{Zufferey2016}
\bibinfo{author}{Zufferey, N.}, \bibinfo{author}{Cho, B.Y.},
  \bibinfo{author}{Glardon, R.}, \bibinfo{year}{2016}.
\newblock \bibinfo{title}{{Dynamic multi-trip vehicle routing with unusual
  time-windows for the pick-up of blood samples and delivery of medical
  material}}, in: \bibinfo{booktitle}{Proceedings of 5th the International
  Conference on Operations Research and Enterprise Systems (ICORES 2016)}, pp.
  \bibinfo{pages}{366--372}.
\newblock \DOIprefix\doi{10.5220/0005733303660372}.

\end{thebibliography}

\newpage
\begin{appendix}
\section{Parameters \label{app_param}}

\begin{table}[H]
    \centering
\caption{Parameters used in the \ac{ALNS} \cite{Christensen2014AdaptiveThesis}}
\label{param_ALNS}
    \begin{tabular}{p{1.5cm} p{7cm} p{2cm}}
    \toprule
    \multicolumn{3}{c}{\textbf{\ac{ALNS} \cite{Rpke2006AnWindows}}}\\
    \midrule
    \textbf{Notation} & \textbf{Description} & \textbf{Value} \\
    \midrule
    $\zeta$ & cooling factor & 0.9999 \\ 
    $\omega$  & start temperature parameter & 0.015 \\
    $\pi$ & maximum removal percentage & 35\% \\
    $\pi'$ & minimum removal percentage & 5\% \\
    $\eta$ & segment size & 125 \\
    $\gamma$ & resetting parameter & 5000 \\
    $r$ & reaction factor & 0.1 \\
    $\sigma_1$ & new global best reward & 33 \\
    $\sigma_2$ & better current solution reward & 9 \\
    $\sigma_3$ & new solution reward & 13 \\
    \midrule
    \multicolumn{3}{c}{Worst removal\cite{Rpke2006AnWindows}} \\
    \midrule
    $p_{worst}$ & probability factor & 4 \\
    \midrule
    \multicolumn{3}{c}{Shaw removal\cite{Rpke2006AnWindows}} \\ 
    \midrule
    $\varphi $ & distance weight & 9 \\ 
    $\chi $ & time weight & 3 \\ 
    $\psi $ & demand weight & 2 \\ 
    $\omega $ & same route weight & 5 \\ 
    $p_{shaw}$ & probability factor & 4 \\
    \midrule
    \multicolumn{3}{c}{\ac{SISR} \cite{ChristiaensJanSIbS}} \\
    \midrule
    $ \overline{c}$ & average number of removed customers & 10 \\
    $ L_{\max}$ & maximum cardinality of removed strings & 10 \\
    $ \alpha$ & probability of not increasing the preserved string size & 0.01 \\
    $\beta$ & blink rate & 0 \\
    \bottomrule
    \end{tabular}
\end{table}

\begin{table}[H]
    \centering
\caption{Parameters for costs by vehicle type}
\label{param_costs_veh_type}
    \begin{tabular}{p{1.5cm} p{2.2cm} p{4cm} p{3cm}}
    \toprule
    \multicolumn{4}{c}{EVs}\\
    \midrule
    \textbf{Notation} & \textbf{Description} & \textbf{Value} & \textbf{Source} \\
    \midrule
    $f_k$ & usage cost in USD & 0 &  \\ 
    $c_k$ &  energy cost in USD per kWh  & 0.1973 for EV, 0.2021 for ICEV & Avg. diesel price is 1.34 USD per liter \\  
    $M_k$ & maintenance cost in USD per km & 0.080837 for electric vans, 0.017015 for electric cargo bikes, and 0.115481 for \acp{ICEV}  & \citet{Mitropoulos2017TotalVehicle} \\
    \bottomrule
    \end{tabular}
\end{table}

\begin{table}[H]
    \centering
\caption{Parameters for efficiency by vehicle type}
\label{param_power_veh_type}
    \begin{tabular}{p{1.5cm} p{4.5cm} p{1.8cm} p{3cm}}
    \toprule
    \multicolumn{4}{c}{\textbf{\acp{EV}}}\\
    \midrule
    \textbf{Notation} & \textbf{Description} & \textbf{Value} & \textbf{Source} \\
    \midrule
    $\phi^d$, $\phi^r$ & electric engine efficiency during discharge and recuperation & 1.184692, 0.846055 & \citet{Goeke2015RoutingVehicles} \\ 
    $\varphi^d $, $\varphi^r $ & coefficient to account for external factors and \ac{SOC} for discharge and recuperation & 1.112434, 0.928465 & \citet{Goeke2015RoutingVehicles} \\
    \toprule
    \multicolumn{4}{c}{\textbf{\acp{ICEV}}}\\
    \midrule
    \textbf{Notation} & \textbf{Description} & \textbf{Value} & \textbf{Source} \\
    \midrule
    $k$ & engine friction factor& 0.2 kJ / rev / l & \citet{Demir2011ATransportation} \\
    $N'$ & engine speed & 33 rev/s & \citet{Demir2011ATransportation} \\
    $D$ &engine displacement & 1.6 liters & \citet{RenaultDanmark2018NyDanmark} \\
    $\eta'$ & efficiency parameter for diesel engines & 0.9 & \citet{Demir2011ATransportation} \\
    $\eta_{tf}$ & drive train efficiency & 0.4 & \citet{Demir2011ATransportation} \\
          \bottomrule
    \end{tabular}
\end{table}

\begin{table}[H]
    \centering
\caption{Parameters for $P_M$ in the energy consumption model}
\label{param_power_Pm}
    \begin{tabular}{p{1.5cm} p{4cm} p{2cm} p{3cm}}
    \toprule
    \multicolumn{4}{c}{$\bm{P_M}$}\\
    \midrule
    \textbf{Notation} & \textbf{Description} & \textbf{Value} & \textbf{Source} \\
    \midrule
        $c_r$ & coefficient of rolling friction & 0.01 & \citet{Demir2011ATransportation} \\ 
      $\alpha$  & gradient of the road & 0\degree & Case study specific. \\
      $\rho_a$ & density of air& $1.2041$ kg/m$^3$ & \citet{Demir2011ATransportation} \\
     $A_f$  & frontal surface area of the vehicle   & $4.06$ m$^2$ & \citet{Kuhlwein2016DRIVING2025} \\
     $c_d$ &   drag coefficient&   $0.34$  & \citet{Kuhlwein2016DRIVING2025} \\
      $m_k(q_j)$ &   total mass due to vehicle in kg $k$ with mass $q_j$ reaching node $j$  & - & Case study specific.  \\
    $g$ &  gravity acceleration  & $9.81$ m/s$^2$  & \citet{Demir2011ATransportation}  \\
      $v_k$ & average speed of vehicle $k$ in m/s &   $12$ m/s & - \\
      $a_k$  & acceleration of vehicle $k$ in m/s$^2$ & 0.01 m/s $^2$ & \citet{Demir2011ATransportation} \\
             \bottomrule
    \end{tabular}
\end{table}
      
\begin{table}[H]
\centering
\small
\caption{Parameters for $P_T$ in the energy consumption model}
\label{param_power_Pt}
\begin{tabular}{p{1.5cm} p{4cm} p{2.5cm} p{3.6cm}}
\toprule
         \multicolumn{4}{c}{$\bm{P_T}$}\\
    \midrule
    \textbf{Notation} & \textbf{Description} & \textbf{Value} & \textbf{Source} \\
    \midrule
      $A$, $B$, $C$    &  constant coefficients from ASHRAE Clear Day Solar Flux Model & 1.088, 0.205, 0.134 & \citet{Ramsey2018SolarRadiation} \\
      $\theta_z$  &  zenith angle (the angle of solar light with respect to the vertical) & 35$^{o}$ & \citet{EUNADICS-AV2018SolarSZA} \\
      $\theta_{i}$ & angle of incidence (with respect to the normal to surface $i$) & 55$^{o}$ & \citet{EUNADICS-AV2018SolarSZA} \\ 
      $\rho_g$ & reflectivity of ground & $0.2$ & \citet{Fayazbakhsh2013ComprehensiveMethod} \\
      $A_i$ & area of surface $i$  & $A_{\text{windshield}}= 1$ m$^2$, $A_{\text{window}} = 0.6$  m$^2$, $A_{\text{rear}} = 0.8$  m$^2$, $A_{\text{roof}} = 1$  m$^2$, $\sum_i A_i = 4$ m$^2$ & \citet{Fayazbakhsh2013ComprehensiveMethod} \\
      $\tau_i$ & emissivity of surface $i$ & $ 0.8$ & \citet{Fayazbakhsh2013ComprehensiveMethod} \\
      $U_i$ & heat transfer coefficient of air with respect to surface $i$  & $0.003$ kW/m$^2 \cdot$K & \citet{Fayazbakhsh2013ComprehensiveMethod} \\
      $C_P$ & specific heat of air at constant pressure &  $1.005$ kJ/kg$\cdot$K & \citet{Fayazbakhsh2013ComprehensiveMethod} \\
      $\eta_{H/C}$ & thermal efficiency of heaters or coolers & $2.5$ for both & \citet{Liu2018ExploringConsumption} \\ 
      $\dot{m}$ & mass of air flow rate from the heater/cooler & $0.3185$ kg/s & \citet{Song2015AirVehicle} \\
      $H_P$ &  human heat production rate & $58.2$ W/m$^2$ & \citet{Fayazbakhsh2013ComprehensiveMethod} \\ 
      $n_p$ & number of passengers & $1$ & - \\ 
      $A_{Du}$ & DuBois area (the surface area of skin of an adult) & $1.8$ m$^2$ & \citet{Fayazbakhsh2013ComprehensiveMethod} \\ 
      $T_d$ & desired temperature & 20\degree C & - \\
      $T_o$ & outside temperature & - & See \Cref{cases} \\
      \bottomrule
    \end{tabular}
\end{table}

\subsection{Parameters for the recharging functions from  \citet{Montoya2017}}

\begin{align*}
\text{30kWh Standard \ac{SOC} in kWh} = \begin{cases}
0.1828 t & 0 \leq t \leq 139.5 \text{ min}
\\ 0.08889 t + 13.1 & 139.5  \leq 173.25   \text{ min} 
\\ 0.02778 t + 23.6875 & 173.25 \leq t \leq 227.25  \text{ min}
\end{cases}
\\ \text{30kWh Fast \ac{SOC} in kWh} = \begin{cases}
0.83089 t & 0 \leq t \leq 30.69  \text{ min}
\\ 0.40404 t + 13.1 & 30.69  \leq 38.115   \text{ min} 
\\ 0.12626 t + 23.6875 & 38.115 \leq t \leq 49.995  \text{ min}
\end{cases}
\end{align*}

\clearpage 
\section{Clarification for constraint (\ref{NoFvisitICEV})} \label{App:ConstraintExplanation}
\begin{figure}[ht!]
        \centering
        \begin{subfigure}[b]{0.45\textwidth}
                \centering
                \includegraphics[width=\textwidth]{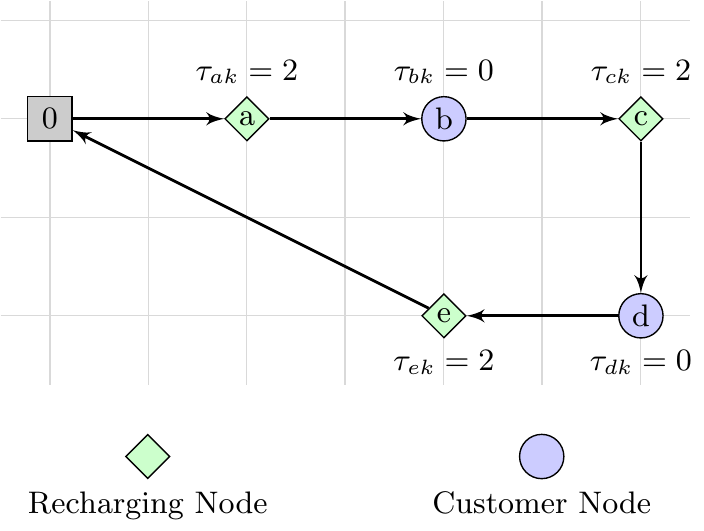}
                \caption{\small{Time can be reset to 0}}
                \label{fig:resetTime}
        \end{subfigure}%
		\qquad
        \begin{subfigure}[b]{0.45\textwidth}
                \centering
                \includegraphics[width=\textwidth]{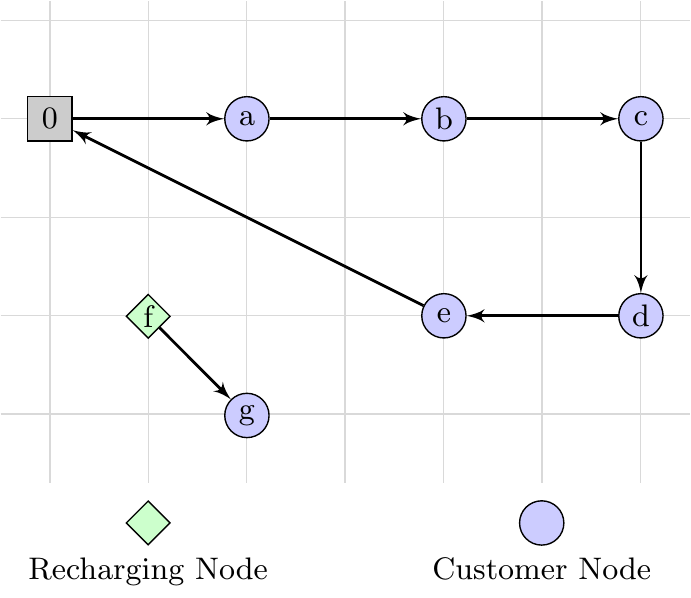}
                \caption{\small{Flow constraints can be circumvented}}
                \label{fig:noFlow}
        \end{subfigure}
        \caption{\small{Examples of undesired solution, feasible in the model of \citet{Goeke2015RoutingVehicles}}}
        \label{fig:UndesiredSolutions}
\end{figure}

In the model of \citet{Goeke2015RoutingVehicles}, there are no constraints enforcing that \acp{ICEV} cannot travel to or from recharging stations. In \Cref{fig:UndesiredSolutions} two feasible solutions for an \ac{ICEV}, $k$ is shown. These two examples show how omitting constraint \eqref{NoFvisitICEV} can lead to undesired solutions.

\Cref{fig:resetTime} shows an example where the time variable can be reset to 0. Constraint \eqref{time} is not posted for arcs leaving an recharging station, and therefore the only thing constraining the $\tau$ variable at the customer nodes is the non-negativity constraint \eqref{nonnegative}. This is obviously a problem as we are working with time-windows.
In \Cref{fig:resetTime} the route alternates between a recharging station and a customer node. That is to satisfy the flow conservation constraint \eqref{flowICEV}. The left hand side of \eqref{flowICEV} does not account for arcs with an end node being a recharging station. Therefore, for all nodes in \Cref{flowICEV} both summation of the left hand side will be equal to $0$.

\Cref{fig:noFlow} shows an example where the flow of an \ac{ICEV} can consist of multiple weakly connected components. Here one component is an actual route starting at the depot and visiting customers $a,b,c,d$, and $e$ before returning to the depot. Another component is the subgraph composed of the nodes $f$ and $g$ in addition to the arc $(f,g)$. As previously described, the flow conservation constraint \eqref{flowICEV} does not account for the arc $(f,g)$, as it originates at a recharging node. Furthermore, the flow conservation constraint is not posted for recharging nodes, making the overall flow feasible for an ICEV.

To disallow these solutions we have included constraint \eqref{NoFvisitICEV} in the model.

\newpage
\section{Instance generation for the Region H case} \label{blood_inst_gen}
\begin{enumerate}
\item Number of customers: Based on the data from Region H, we created a truncated normal distribution whose output is rounded to the closest integer. It has the following parameters:
\begin{itemize}
\item	Mean = 107.8946,
\item Standard deviation = 26.63986, and
\item Truncation (lowest: 26.5, highest: 145.5).
\end{itemize}

\item Taxi cost: The average cost of taxi trips is estimated to be 30 USD since a taxi trip in Copenhagen costs around 30 USD for 8km, which is twice the distance from Bispebjerg Hospital to Frederiksberg Hospital, based on information from the DanTaxi website 
\\(\url{https://dantaxi4x48.dk/hvad-koster-en-taxa/}).

\item Temperature: We calibrated four triangular distributions, one for each season, based on data in degrees Celsius from the Danish Meteorological Institute (\url{dmi.dk}). Each season has a 25\% probability of being selected.
\begin{table}[ht!]
\centering
\caption{Description of seasonal triangular distributions of temperature}
\label{temp_dist_blood}
\begin{tabular}{lrrr}
\toprule
	&\textbf{Minimum}	& \textbf{Maximum}	& \textbf{Mean} \\
	\midrule
Summer &	13.4 &	26.9	& 19.5
\\ Autumn &	1.2 &	21.8 &	10.7
\\ Winter &	-7.3 &	8.8 &	2.3
\\ Spring &	-4	& 24.1 &	8.2 \\
\bottomrule
\end{tabular}
\end{table}
\item List of customers' locations: Each clinic has a certain probability of being selected, which was calculated from the historical data provided by Region H. Addresses of the clinics  were provided by Region H (also available publicly) and were geocoded using the tool provided at the website \url{https://geocode.localfocus.nl/}.

\item Customer time windows and service times: We obtained the time windows and service times (of 3 minutes per clinic) for various clinics based on data provided by Region H. 

\item Demands: Based on an online survey conducted by us, we calculated the average number of vials handled by a doctor depending on the day of the week and on whether it is summer (lower due to vacation) or the rest of the year, so in total we calculated 10 averages (summer Monday, rest of the year (ROTY) Monday, summer Tuesday, etc.). We rounded those averages to the closest integer and multiplied them by the number of doctors in each clinic to obtain the demand per clinic. The days of the week were selected with a 20\% probability each while summer was selected with a 25\% probability.
\begin{table}[ht!]
\centering\caption{Demand quantity per doctor per day by season (in number of vials)}
\begin{tabular}{cccccc}
\toprule
\multirow{2}{*}{Summer} & Monday &Tuesday & Wednesday& Thursday &Friday \\
&6 &6 &	4 &	6 & 5 \\
\midrule
\multirow{2}{*}{ROTY} &Monday &Tuesday & Wednesday &Thursday & Friday \\
& 9 & 8 & 6 & 7 & 6 \\
\bottomrule
\end{tabular}
\end{table}
\item Travel times: Free flow travel times between origin-destination pairs are obtained from OSRM \cite{luxen-vetter-2011} and scaled up by zone-dependent congestion factors obtained from the Danish National Transport Model \cite{Rich2016TheResults}. 
\end{enumerate}
\section{Instance generation for the MTH case}
\label{mth_inst_gen}

\begin{enumerate}
\item Number of customers: Based on data provided by MT Højgaard, we created a truncated normal distribution whose output is rounded to the closest integer. It has the following parameters:
\begin{itemize}
\item	Mean = 423.3636,
\item Standard deviation = 34.2547, and
\item Truncation (lowest: 348, highest: 476).
\end{itemize}

\item Taxi cost: The cost of a taxi in Greater Copenhagen for an average distance of 10 km is around 30 USD, according to Dantaxi (\url{https://dantaxi4x48.dk/hvad-koster-en-taxa/}).

\item Temperature: We calibrated four triangular distributions, one for each season, based on data in degrees Celsius from the Danish Meteorological Institute (\url{dmi.dk}). Each season has a 25\% probability of being selected.
\begin{table}[ht!]
\centering
\caption{Description of seasonal triangular distributions of temperature}
\label{temp_dist_blood1}
\begin{tabular}{lrrr}
\toprule
	&\textbf{Minimum}	& \textbf{Maximum}	& \textbf{Mean} \\
	\midrule
Summer &	13.4 &	26.9	& 19.5
\\ Autumn &	1.2 &	21.8 &	10.7
\\ Winter &	-7.3 &	8.8 &	2.3
\\ Spring &	-4	& 24.1 &	8.2 \\
\bottomrule
\end{tabular}
\end{table}

\item List of customers' locations: Each customer has a certain probability of being selected, which was calculated from the historical data provided by MT Højgaard. Anonymous customer addresses were provided by the company and were geocoded using the tool provided at the website \url{https://geocode.localfocus.nl/}.

\item Service times and time windows: A normal probability distribution for the service time per customer visit was calculated from data provided by MT Højgaard. This distribution was applied to the aforementioned list of customers, and therefore a static service time per customer was obtained. Time windows were estimated in the following way:
\begin{itemize}
\item If the customer visit lasts less than 2 hours, it is assigned either to the morning (first 5 hours of the shift, from 6 to 11 am) or the afternoon (last 5 hours of the shift, 6 am to 4 pm).
\item If the customer visit lasts more than 2 hours, it can be performed anytime within the 10-hour shift (6 am to 4 pm). Note that the shifts are 10 hours long because they include transportation to and from home as well as breaks, which are modeled as customer visits.
\end{itemize}

\item Demands: Since the nature of this problem deals with the routing of technicians, we assume they always carry the same equipment in their vehicles and therefore we consider the demand to be null for all customer visits.

\item Compatibility matrix: We followed a nested approach for the computation of the compatibility matrix. Following guidelines from Lindpro, we consider four driver types/technician skill sets and four different task types depending of the level of the skills required. 
\begin{table}[h!]
\centering
\caption{Compatibility between technician skill sets (percentage over total number of drivers) and task types (percentage over total number of tasks)}
\label{comp_matrix}
\begin{tabular}{lcccc}
\toprule
	&\textbf{Level 1 (40\%)}	& \textbf{Level 2 (30\%)}	&\textbf{Level 3 (20\%)} & \textbf{Level 4 (10\%)} \\
	\midrule
Classified Technician (10\%)   &  &   &X   &X
\\ Technician (20\%)    &  &X   &X   &
\\ Electrician (30\%)   &  &X   &   &
\\ Handyman (40\%) &X  &   &   & \\
\bottomrule
\end{tabular}
\end{table}
Technicians trained to a certain level are compatible with tasks on their skill level and the immediately inferior one with the exception of handymen, who have their own task category. \Cref{comp_matrix} illustrates our approach for the calculation of the compatibility matrix with the different skill sets and task types in our model. Note that the higher the task type number, the higher the level of skills required.

\item Driver home locations: Anonymous addresses corresponding to the driver's home locations were provided by MT Højgaard and were geocoded using the tool provided at the website \url{https://geocode.localfocus.nl/}.

\item \ac{EV} eligibility: We consider a vehicle to be EV-eligible if the home location of its driver is less than 1 km away from the nearest EV charging station. Both the networks of E.ON and Clever in the urban area of Copenhagen are included in this consideration. This approach resulted in 44 EV-eligible vehicles. 

\item Travel times: Free flow travel times between origin-destination pairs are obtained from OSRM \cite{luxen-vetter-2011}, which are then suitably scaled up by the corresponding zone-dependent congestion factors obtained from the Danish National Transport Model \cite{Rich2016TheResults}. 

\end{enumerate}

\end{appendix}

\end{document}